\documentclass[12pt]{article}

\usepackage{fullpage}

\usepackage{amsfonts}
\usepackage{graphicx}
\usepackage{mathrsfs}
\usepackage{amsmath}
\usepackage{amssymb}
\usepackage{float}
\usepackage{color}
\usepackage{booktabs}
\usepackage[authoryear]{natbib}
\usepackage{verbatim}

\renewcommand{\baselinestretch}{1.7}
\usepackage{geometry}
\newtheorem{pro}{Proposition}

\newtheorem{thm}{Theorem}
\newtheorem{assu}{Assumption}
\newtheorem{lemma}{Lemma}
\def\D{{\mathcal D}}
\def\bH{\mathbb{H}}
\numberwithin{equation}{section}
\newtheorem{coro}{Corollary}
\def\cp{\mathop{\rightarrow}\limits^{p}}

\def\remark#1{\noindent{\bf Remark #1\ }}

\def\proof {{\noindent\bf Proof.}\quad}
\newcommand{\tc}[1]{\multicolumn{1}{c}{#1}}

\def\var{\mathrm {Var}}

\def\E{\mathbb {E}}

\def\D{{\mathcal D}}
\def\A{{\mathcal A}}

\def\ca{{\rm Card}}

\def\tW{\tilde{W}}
\def\wht{{L}}

\def\I{\mathcal{I}}
\def\bI{\mathbb{I}}
\def\bH{\mathbb{H}}

\newcommand{\bm}{\boldsymbol}
\def\FDR{\mathrm {FDR}}
\def\FDP{\mathrm {FDP}}
\def\var{\mathrm {Var}}
\def\sgn{\mathrm {sgn}}

\date{}

\begin{document}

\renewcommand{\baselinestretch}{1.0}

\date{}
\title
{\bf A New Procedure for Controlling False Discovery Rate in Large-Scale $t$-tests}
\author{Changliang Zou$^a$, Haojie Ren$^b$, Xu Guo$^c$ and Runze Li$^b$\\
$^a$School of Statistics and Data Science, Nankai University\\
 Tianjin, P.R. China\\
$^b$Department of Statistics, The Pennsylvania State University, \\
University Park, PA 16802-2111, USA\\
$^c$School of Statistics, Beijing Normal University\\
 Beijing, P.R. China
}

\maketitle

\begin{abstract}
This paper is concerned with false discovery rate (FDR) control in large-scale multiple
testing problems. We first propose a new data-driven testing
procedure for controlling the FDR in large-scale $t$-tests for one-sample mean
problem.
The proposed procedure
achieves exact FDR control in finite sample settings when the populations are
symmetric no matter the number of tests or sample sizes. Comparing
with the existing bootstrap method for FDR control, the proposed procedure is computationally efficient.
%
We show that the proposed method can control the FDR asymptotically for asymmetric
populations even when the test statistics are not independent. We further show that
the proposed procedure with a simple correction is as accurate as the bootstrap method to
the second-order degree, and could be much more effective than the existing normal calibration.
We extend the proposed procedure to two-sample mean problem.
Empirical results show that the proposed procedures have better FDR control than existing
ones when the proportion of true alternative hypotheses is not too low, while
maintaining reasonably good detection ability.
\noindent
\\[6pt]
\noindent{\bf Keywords}: Data splitting; Large-deviation probability;  Multiple comparisons; Product of two normal variables; Skewness; Symmetry
\end{abstract}

\section{Introduction}


Many multiple testing problems are closely related to one-sample mean problem. Let
$(X_{i1},$ $\ldots,X_{ip})^{\top}$, $1\leq i\leq n_t$ be independent and identically distributed (i.i.d.)
samples from $X=(X_1,\ldots,X_p)^{\top}$ with mean $\mu=(\mu_1,\ldots,\mu_p)^\top$.
Of interest is to test $H_0:\mu=0$ versus $H_1: \mu\ne 0$. This leads to
consider a multiple
testing problem on the mean values
\[
 \mathbb{H}_{0j}: \mu_j=0\ \ \mbox{v.s.} \  \ \mathbb{H}_{1j}: \mu_j\neq 0, \ 1\leq j\leq p.
\]
A standard procedure for false discovery rate (FDR) control
is to apply the Benjamini and Hochberg (BH) method \citep{benjamini1995controlling}
to the $t$ statistics with the standard normal or Student's $t$ calibration. That is,
let $\bar{X}_j=n_t^{-1}\sum_{i=1}^{n_t}X_{ij}$ and $s_j^2=(n_t-1)^{-1}\sum_{i=1}^{n_t}(X_{ij}-\bar{X}_j)^2$, and
define $T_j=\sqrt{n_t}\bar{X}_j/s_j$, then the procedure
rejects a hypothesis whenever $|T_j|\geq T$ with a data-dependent threshold
\begin{align}\label{bht}
T=\inf\left\{t>0: \frac{2p\{1-\Phi(t)\}}{\#\{j:|T_j|\geq t\}}\leq \alpha\right\},
\end{align}
for a desired FDR level $\alpha$ \citep{storey2004strong},
 where  $\Phi(\cdot)$ is the cumulative distribution function of $N(0,1)$. It has been revealed that the accuracy of
this control procedure heavily depends on the skewness of $X_j$'s and the diverging rate
of $p$ relatively to $n$, since $\Phi(t)$ is only an approximation to the distribution of
$T_j$.

Many studies have investigated the performance of large-scale $t$ tests. \cite{efron2004large} observed that the null distribution
choices substantially affect the simultaneous inference procedure in a microarray
analysis.
\cite{delaigle2011robustness} conducted a careful study of
moderate and large deviations of the $t$ statistic which is indispensable to
understanding its robustness and drawback for analyzing high dimensional data. Under a
condition of non-sparse signals, \cite{cao2011simultaneous} proved the robustness
of Student's $t$ test statistics and $N(0, 1)$ calibration in the control of FDR. \cite{liu2014phase} gave a systematic analysis on the asymptotic conditions with which the large-scale $t$ testing procedure is able to have FDR control.

Bootstrap method is known as a useful way to improve the accuracy of an exact null
distribution approximation and has been demonstrated to be particularly effective for
highly multiple testing problems. See 
\cite{delaigle2011robustness} and
the references therein. In general, the bootstrap is capable of correcting for skewness,
and therefore leads to second-order accuracy. Accordingly, a faster increasing rate of
$p$ could be allowed \citep{liu2014phase} and better FDR 
control would be
achieved when the data populations are skewed. However, multiple testing problems with
tens of thousands or even millions of hypotheses are now commonplace, and practitioners
may be reluctant to use a bootstrap method in such situations, and therefore a rapid
testing procedure is highly desirable.

In this paper, we propose a new data-driven selection procedure
controlling the FDR. The method entails constructing $p$ new test statistics
with marginal symmetry property, using the empirical distribution of
the negative statistics to approximate that of the
positive ones, and searching for the threshold with a formula similar to
(\ref{bht}). The proposed procedure is computationally efficient since it only uses a
one-time split of the data and calculation of the product of two $t$ statistics
obtained from two splits. We study theoretical properties  of the proposed procedure. We show that (a)
the proposed method achieves exact FDR control even in finite sample settings
when $X_j$'s are symmetric and independent; and (b) the proposed method achieves
asymptotical FDR control under mild conditions when the populations are
asymmetric and dependent. We further propose a simple refinement of the proposed procedure,
and study the asymptotical property of the refined one. The theoretical property
of the proposed refined one implies that
it is as accurate as the bootstrap method to
the second-order degree in certain situations.
We also investigate extension of the proposed procedure to two-sample mean problem.
Simulation comparisons imply that the proposed method has better FDR control
than existing methods, maintains reasonably good power and has a significant
reduction in computing time and storage.

The rest of this paper is organized as follows. In Section 2, we present
the new procedure and establish its FDR control property.  Some extensions are given in
Section 3. Numerical studies are conducted in Section 4. Section 5 concludes
the paper, and theoretical proofs are delineated in the Appendix. Some
technical details and additional numerical results are provided in the
Supplementary Material.

Notations. $A_n\approx B_n$ stands for that $A_n/B_n\cp 1$ as $n\to\infty$.
The ``$\gtrsim$"
and ``$\lesssim$" are similarly defined. We denote by $\I_0$ and $\I_1$ the
true null set and alternative set, $p_0=|\I_0|$ and $p_1=|\I_1|$, respectively.

\section{A New FDR Control Procedure and its Theoretical Properties}

We first propose a new FDR control procedure for the one-sample mean problem,
and then establish the theoretical properties of the proposed procedure.

\subsection{A new FDR control procedure}

Without loss of generality, assume that the sample size is an even integer $n_t=2n$.
We randomly split the data into two disjoint groups $\D_{1}$ and $\D_{2}$ of equal
size $n$. The $t$ test statistics for the $j$th variable on $\D_{1}$ and $\D_{2}$ are denoted as $T_{1j}$ and $T_{2j}$, respectively.
Define
\[
W_j=T_{1j}T_{2j}.
\]
Clearly, $W_j$ is likely to be large for most of the signals regardless of the sign of $\mu_j,j\in\I_1$, and small for most of the null
variables. Observing that $W_j$ is, at least asymptotically, symmetric with mean zero for $j\in\I_0$ due to the central limit theorem and the independence between $\D_{1}$ and $\D_{2}$, we can choose a threshold $L>0$ by setting
\begin{align}\label{th}
L=\inf\left\{t>0:\frac{\#\{j: W_j\leq -t\}}{\#\{j:W_j\geq t\}\vee 1}\leq \alpha\right\},
\end{align}
and reject the $\bH_{0j}$ if $W_j\geq L$, where $\alpha$ is the target FDR
level. If the set is empty, we simply set $L=+\infty$. The fraction in
(\ref{th}) is an estimate of the false discovery proportion (FDP) since the set
$\{j: W_j\leq -t,j\in\I_1\}$ is often very small (if the signal is not too
weak) and thus $\#\{j: W_j\leq -t\}$ is a good approximation to $\#\{j: W_j\leq
-t, j\in\I_0\}$.

\begin{figure}[ht]\centering
\includegraphics[width=1\textwidth]{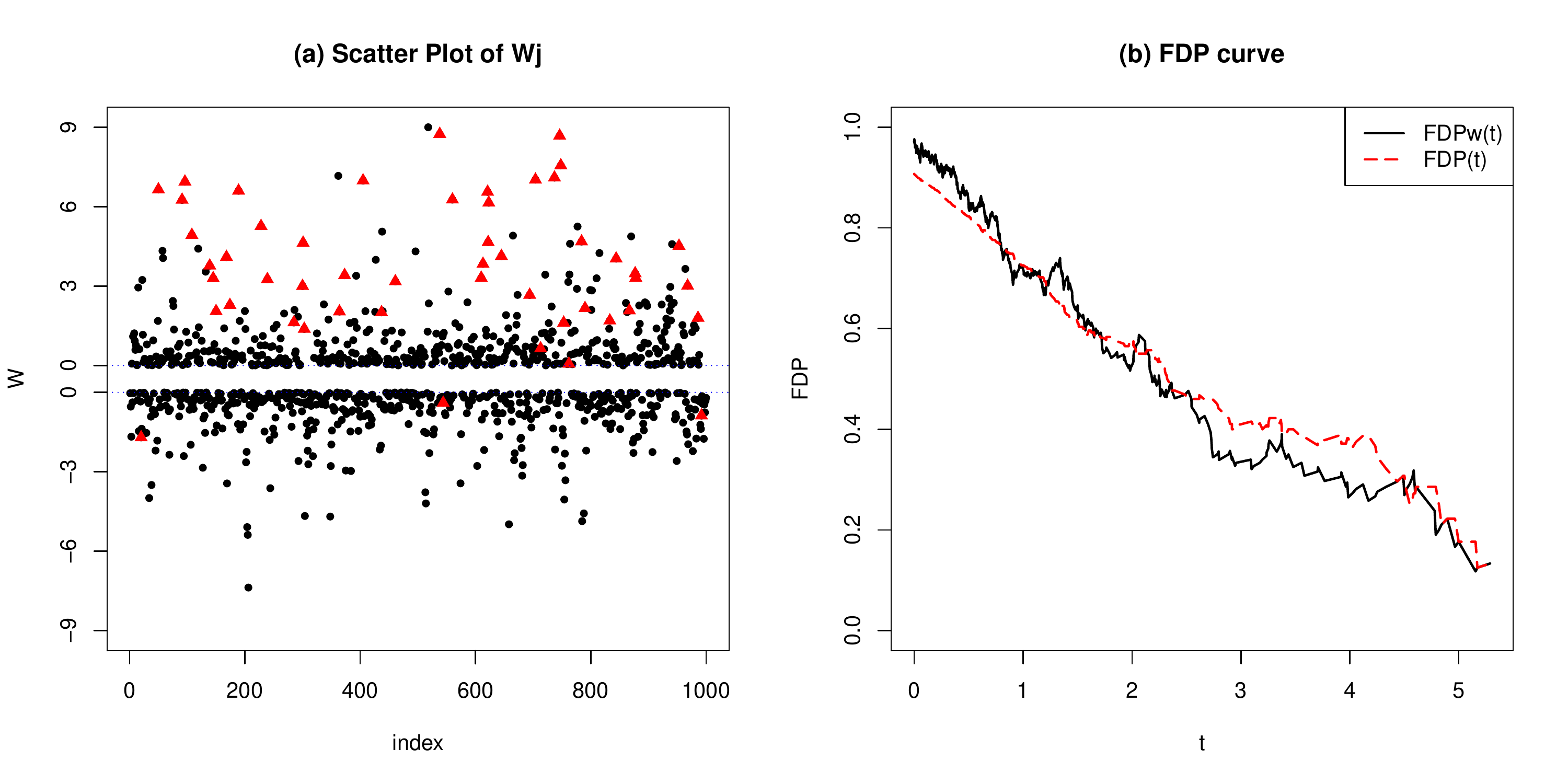}\vspace{-0.2cm}
\caption {\small \it (a): Scatter plot of the $W_j$'s with red triangles and black dots denoting true signals and nulls respectively; (b): The corresponding estimate of FDP curve (against $t$) along with the true FDP. In this case, $p=1000$, $n_t=100$ and $|\I_1|=50$. We consider a multivariate chi-squared distribution with two degrees of freedoms and an autoregressive structure $(0.5^{|i-j|})$.  We set the signal-to-noise ratio as 2.}\label{Fig:tst}\vspace{-0.2cm}
\end{figure}

As described above, we construct the test statistic $W_j$ with marginal
symmetry property by using sample splitting. Thus we refer this procedure to as {\it Reflection via Sample Splitting} (RESS). The RESS procedure is {\it data-dependent} and does not
involve any unknown quantity related to the populations. This is an important
and desirable feature of the RESS. Figure \ref{Fig:tst} depicts a visual
representation of the RESS procedure. Specifically, Figure \ref{Fig:tst}(a)
depicts the scatter plot of the $W_j$'s with red triangles denoting true
signals. Observe that the true signals are primarily above the x-axis,
indicating $W_j>0$, while the null $W_j$'s (black dots) are roughly
symmetrically distributed across the horizontal lines. Figure \ref{Fig:tst}(b)
depicts the corresponding estimate of FDP (i.e., the fraction in (\ref{th}))
along with the true FDP over $t$. The approximation in this case is very good
as only three true alternatives (i.e., three red triangles) lie below the
horizontal line in Figure~\ref{Fig:tst}(a).

Knockoff framework was introduced by \cite{bc2015} in the high-dimensional
linear regression setting. The knockoff selection procedure operates by
constructing ``knockoff copies" of each of the $p$ covariates (features) with
certain knowledge of the covariates or responses, which are then used as a
control group to ensure that the model selection algorithm is not choosing too
many irrelevant covariates. The signs of test statistics in the knockoff need
to satisfy (or roughly) joint exchangeability so that the corresponding
knockoff can yield accurate FDR control in finite sample setting. Refer to
\cite{candes2018panning} for more discussions. The proposed threshold $L$ in (\ref{th})
shares a similar spirit to the knockoff, but they are distinguished in that the
RESS procedure does not require any prior information about the distribution
of $X=(X_1,\ldots,X_p)^{\top}$. This is especially important since it is
difficult to estimate the distribution of $X$ when $p$ is very large. We employ
the sample-splitting strategy to achieve a marginal symmetry property. It turns
out that the FDR can be controlled reasonably well due to the marginal
symmetry of $W_j$'s. The theoretical findings on the FDR control under certain
dependence structures such as positive regression dependence on subset or weak
dependence at a marginal level \citep{benjamini2001control,storey2004strong}
shed light on the validity of the RESS procedure. Detailed
analysis will be given in Section 2.2.

\begin{figure}[ht]\centering
\includegraphics[width=0.9\textwidth]{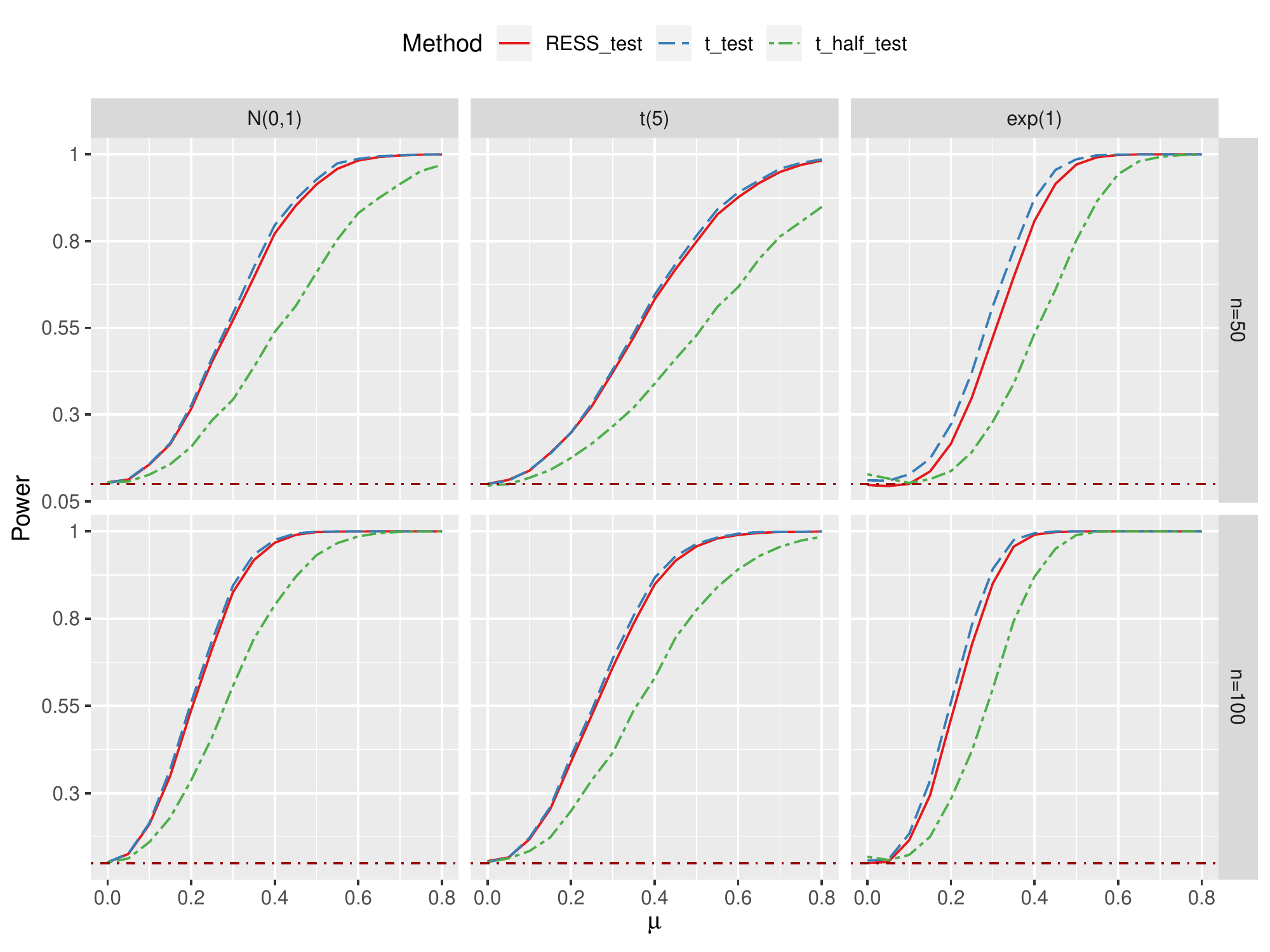}\vspace{-0.2cm}
\caption {\small \it Power comparison of the tests based on $W_j$ (solid line), $T_j$ (dashed line) and $T^2_{1j}$ (dashed-dotted line) when $n_t=50$ or $100$. The type I error is fixed as 0.05 and
the error distributions considered include $N(0,1)$, the Student's $t$ distribution with five degrees of freedom $t(5)$ and the exponential distribution.}\label{Fig:tst1}\vspace{-0.2cm}
\end{figure}

At a first glance, the test statistic $W_j$ may result in much information loss due to the sample-splitting. In fact, benefiting from the joint use of two independent $t$ statistics,
the relative power loss of $W_j$ with respect to $T_j$ is quite mild. By the inequality $\Pr^2(T_{1j}>\sqrt{t})+\Pr^2(T_{1j}<-\sqrt{t})\leq \Pr(T_{1j}T_{2j}>t)\leq 2\Pr(T_{1j}>\sqrt{t})+2\Pr(T_{1j}<-\sqrt{t})$,  the power ratio of the tests based on $W_j$ and $T_j$ can be easily bounded as
\begin{align*}
&\frac{\left\{1-\Phi(z_{\alpha/4}-\sqrt{n}\mu_j/\sigma_j)\right\}^2+\left\{1-\Phi(z_{\alpha/4}+\sqrt{n}\mu_j/\sigma_j)\right\}^2}{2-\Phi(z_{\alpha/2}-\sqrt{2n}\mu_j/\sigma_j)-\Phi(z_{\alpha/2}+\sqrt{2n}\mu_j/\sigma_j)} \lesssim\frac{\Pr_{\bH_{1j}}(W_j>W_{\alpha})}{\Pr_{\bH_{1j}}(|T_j|>z_{\alpha/2})}\\
&\lesssim \frac{2\left\{2-\Phi(z_{\sqrt{\alpha/2}}-\sqrt{n}\mu_j/\sigma_j)-\Phi(z_{\sqrt{\alpha/2}}+\sqrt{n}\mu_j/\sigma_j)\right\}}{2-\Phi(z_{\alpha/2}-\sqrt{2n}\mu_j/\sigma_j)-\Phi(z_{\alpha/2}+\sqrt{2n}\mu_j/\sigma_j)},
\end{align*}
where $\sigma_j$ is the standard deviation of $X_j$, and $W_{\alpha}$ and $z_{\alpha}$ are the upper $\alpha$ quantiles of the distributions of $W_j$ and $N(0,1)$, respectively. Further when $\sqrt n\mu_j\to \infty$, both of the test statistics $W_j$ and $T_j$ have asymptotic power 1. For better understanding, the power curves (with size corrected) of the two tests are presented in Figure \ref{Fig:tst1} for some commonly used settings. We can see that though $W_j$ is always inferior to $T_j$ as we can expect, the disadvantage is not very significant and also tends to be less pronounced when $n$ increases. This power sacrifice of the RESS in turn brings us much better error rate control as we shall show in the next subsection. On the other hand, compared with the test statistic $T_{1j}^2$ based on only group $\D_1$, the proposed test statistic $W_j=T_{1j}T_{2j}$ has smaller variance because $\var(W_j)\approx 1<2\approx \var(T^2_{1j})$. Since the null distribution of $W_j$ is symmetric, the upper quantiles of the distributions of $W_j$ would be much smaller than those of $T^2_{1j}$. As a result, $W_j$ is more powerful than $T_{1j}^2$.

\remark 1 It is noteworthy that the joint use of $W_j$ and the threshold $L$
distinguishes our RESS procedure from the methods given by
\cite{wasserman2009high} and \cite{meinshausen2009p} which used the
sample-splitting scheme to conduct variable selection with error rate control.
They used the first split to reduce the number of variables to a manageable
size and then applied FDR control methods on the remaining variables using the
data from the second split. The normal calibration is usually required to
obtain $p$-values. In contrast, in the RESS procedure, the data from both two
splits are used to compute the statistics and the empirical distribution is in
place of asymptotic distributions.

\subsection{Theoretical results}

We firstly investigate the FDR control of the proposed RESS procedure when
$X_1,\ldots,X_p$ are independent each other, and then extend the results to the
dependent case under stronger conditions. A simply yet effective refined
procedure with better accuracy in FDR control will be further developed after
the convergence rate of the FDP of the RESS is investigated.

\subsubsection{Independent case}

A preliminary result of this paper is that the proposed procedure controls a quantity nearly equal
to the FDR when the populations are symmetric.
\begin{pro}\label{pro1}
Assume $X_1,\ldots,X_p$ are symmetrically distributed and independent each other. For any $\alpha\in(0,1)$ and $n_t\geq 4$, the RESS method satisfies
\[
 \E\left[\frac{\#\{j: W_j\geq L,j \in\I_0\}}{\#\{j: W_j\geq L\}+\alpha^{-1}}\right]\leq \alpha.
\]
\end{pro}
The term bounded by this proposition is very close to the FDR in settings
where $\alpha^{-1}$ is dominated by $\#\{j: W_j\geq L\}$. Following \cite{bc2015}, if it is preferable to control the FDR exactly, we may adjust the threshold by
\[
L_{+}=\inf\left\{t>0:\frac{1+\#\{j: W_j\leq -t\}}{\#\{j:W_j\geq t\}\vee 1}\leq \alpha\right\},
\]
with which we can show that
\[
\mbox{FDR}_{W}(L_{+})=\E\left[\frac{\#\{j: W_j\geq L_{+},j \in\I_0\}}{\#\{j: W_j\geq L_{+}\}}\right]\leq \alpha.
\]
In what follows, as we mainly focus on the asymptotic FDR control, the results
with $L$ and $L_{+}$ are generally the same. From the proof of this
proposition, we can see that the inequality is due to the fact that
\[
 \#\{j: W_j\leq -t,j \in\I_0\}\leq  \#\{j: W_j\leq -t\}
\]
which would usually be tight because most strong signals yield a positive $W_j$ or at least a not too small value of $W_j$. This implies that it is very likely that the FDR of the RESS will be fairly close to the nominal one unless a large proportion of $\mu_j$'s for $j\in\I_1$ is very weak.

Proposition \ref{pro1} is a direct corollary of the following result in which $X_j$'s are allowed to be asymmetric.
\begin{pro}\label{pro2}
Assume $X_1,\ldots,X_p$ are independent each other and $n_t\geq 4$. For any $\alpha\in(0,1)$, the RESS method satisfies
\[
\E\left[\frac{\#\{j: W_j\geq L,j \in\I_0\}}{\#\{j: W_j\geq L\}+\alpha^{-1}}\right]\leq \min_{\epsilon\geq 0}\left\{\alpha(1+4\epsilon)+\Pr\left(\max_{j\in\I_0}\Delta_j>\epsilon\right)\right\},
\]
where $\Delta_j=\left|\Pr(W_j>0\mid |W_j|)-1/2\right|$.
\end{pro}

We can interpret $\Delta_j$ as measuring the extent to which the symmetry is violated for a specific variable $j$. This result concurs with our intuition that controlling the $\Delta_j$'s is sufficient to ensure control of the FDR for the RESS method. In the most ideal case where $X_1,\ldots,X_p$ are symmetrically distributed,  $\Delta_j= 0$ for all $j\in \I_0$, and we automatically obtain the FDR-control result in Proposition \ref{pro1} since we can take $\epsilon=0$. Under asymmetric scenarios, the $\Delta_j$ can still be expected to be small due to the convergence of $T_{1j}$ and $T_{2j}$ to the normal if $n$ is not too small.  In the next theorems, we will show that under mild conditions $\max_{j\in\I_0}\Delta_j\to 0$ in probability, yielding a meaningful result on FDR control in more realistic settings. The proof of this proposition follows similarly to Theorem 2 in \cite{barber2018robust} which shows that the Model-X knockoff \citep{candes2018panning} selection procedure incurs an inflation of the FDR that is proportional to the errors in estimating the distribution of each feature conditional on the remaining features.

For our asymptotic analysis, we need the following assumptions. Throughout this paper, we assume $p_1\leq \gamma p$ for some $\gamma<1$, which includes the sparse setting $p_1=o(p)$.

\begin{assu}\label{moment}
(Moments) (i) For some constant $K_1>0$, $\max_{1\leq j\leq p}\mathbb{E}(X_j-\mu_j)^4/\sigma_j^4<K_1$; (ii) For some constants $C>0$ and $K_2>0$, $\max_{1\leq j\leq p}\mathbb{E}[\exp\{C(X_j-\mu_j)^2/\sigma_j^2\}]<K_2$.
\end{assu}

\begin{assu}\label{signall}
(Signals) $\beta_p\equiv\ca\left\{j: |\mu_j|/\sigma_j\geq 3\sqrt{\log p/n}\right\}\to \infty.$
\end{assu}

\remark 2 The moment condition in Assumption \ref{moment}-(i) is required in a large deviation result for the Student's $t$ statistics on which our proof heavily hinges.
Assumption \ref{moment}-(ii), which requires exponentially light tails and implies that all moments of $X$ are finite, is stronger than Assumption \ref{moment}-(i). This will be only needed when we want to use the RESS method with correction (see Section 2.2.2). In fact, for the familywise error control with bootstrap calibration, similar conditions are also imposed to achieve better accuracy \citep{FHY2007}. The implication of Assumption \ref{signall} is that $p_1\to\infty$. If the number of true alternatives is fixed as $p\to\infty$, \cite{liu2014phase} have shown that even with the true $p$-values, the BH method is unable to control FDP with a high probability. Thus, we use this condition to rule out such cases. \hfill$\Box$

\begin{thm}\label{thm1}
Assume $X_1,\ldots,X_p$ are independent each other.
Suppose Assumptions \ref{moment}-(i), \ref{signall} and $p=o\{\exp(n^{1/3})\}$ hold.  For any $\alpha\in(0,1)$, the FDP of the RESS method satisfies
\begin{align}
{\rm FDP}_W(L)&\equiv\frac{\#\{j: W_j\geq L,j \in\I_0\}}{\#\{j: W_j\geq L\}\vee 1}\nonumber \\
&\leq\alpha+O_p\left\{\sqrt{\xi_{n,p}/n}+n^{-1}(\log p)^{3}\max_{j\in \I_0}\kappa_j^2+\beta_p^{-1/2}\right\}, \label{fa}
\end{align}
and $\mathop{\lim\sup}_{(n,p)\to\infty}{\rm FDR}\leq \alpha$, where $\kappa_j=\E\{(X_j-\mu_j)^3\}/\sigma_j^3$ and $\xi_{n,p}=\max(\log p,\log n)$.
\end{thm}
The proof of this theorem relies on a nice large-deviation result for $t$-statistics \citep{delaigle2011robustness} but our new statistic $W_j$
makes that the technical details of our theory are not straightforward and cannot be obtained from existing works. Under a finite fourth moment condition, Theorem \ref{thm1} reveals that the FDR of our RESS method can be controlled if $\I_1$ satisfies the technical condition on the alternative set, Assumption \ref{signall}. Roughly speaking, the theorem ensures that the $\max_{j\in\I_0}\Delta_j$ in Proposition \ref{pro2}, which can be bounded by $\max_{j\in\I_0}\sup_{0\leq t\leq {2\log p}}|f_j(t)/f_j(-t)-1|$, is small, where $f_j(\cdot)$ denotes the probability density function of $W_j$.
Note that the inequality in (\ref{fa}) is mainly due to
\begin{align*}
\frac{\sum_{j\in\I_0}\bI\left(W_j\geq L\right)}{\sum_{j}\bI\left(W_j\leq-L\right)}\approx \frac{\sum_{j\in\I_0}\bI\left(W_j\leq -L\right)}{\sum_{j}\bI\left(W_j\leq-L\right)}\leq 1.
\end{align*}
Hence, the FDR control is often quite tight because
\[
\sum_{j\in\I_1}\bI\left(W_j\leq -L\right)/\sum_{j\in\I_0}\bI\left(W_j\leq -L\right)\cp 0
\]
in many situations, such like $(p_1-|\mathcal{M}|)/p_0\to 0$, where $\mathcal{M}=\left\{i: |\mu_i|/\sigma_i\geq \sqrt{(2+c)\log p/n}\right\}$ for any small $c>0$. See a proof given in the Supplementary Material.

It is interesting to further unpack the convergence rate given in this theorem. \cite{liu2014phase} has shown that the convergence rate of the bootstrap calibration is
\begin{align*}
 {\rm FDP}_{\rm Bootstrap}&\lesssim\alpha+O_p\left\{\sqrt{\log p/n}+n^{-1}(\log p)^{2}\right\}.
\end{align*}
The ${\rm FDR}_{W}$ indicates that our ``raw" RESS method is inferior to the bootstrap calibration, especially when $p$ is very large (so that the term of order $n^{-1}$ has non-ignorable effect). Actually, the term $n^{-1}(\log p)^3$ can be eliminated by a simple correction as discussed below.

\subsubsection{A refined procedure}
By more carefully examining the FDP of the proposed RESS method, we can show
that for any $0\leq t\leq 2\log p$, we have
\begin{align}\label{FDPW}
\widehat{\rm FDP}_W(t)&\equiv\frac{\#\{j: W_j\leq -L\}}{\#\{j: W_j\geq L\}\vee 1}\nonumber\\
&\approx {\rm FDP}_W(t)\left(1-\frac{2t^3\bar{\kappa}}{9n}\right)+O_p\left(\sqrt{\xi_{n,p}/n}+\beta_p^{-1/2}\right),
\end{align}
where $\bar{\kappa}=p_0^{-1}\sum_{j\in\I_0}\kappa_j^2$. Say, we are able to express the term of order $n^{-1}t^3$ to a more accurate way, which benefits from utilizing the empirical distribution $p^{-1}\sum_{j}\bI(W_j\leq-t)$ to approximate $p_0^{-1}\sum_{j\in\I_0}\bI(W_j\geq t)$ and ``surprisingly" eliminate the terms of order $n^{-1/2}(\log p)^{3/2}$. Clearly, the $\widehat{\rm FDP}_W(t)$ is an underestimate of the true FDP, and in turn yields an inflation of the FDR.

This motives us to consider the test statistic
$\tilde{W}_j=n\bar{X}_{1j}\bar{X}_{2j}/s_{1j}^2=T_{1j}\tilde{T}_{2j}$, where
$\tilde{T}_{2j}=\sqrt{n} \bar{X}_{2j}/s_{1j}$. As shown in the Appendix, the
FDP of using $\tilde{W}_j$ satisfies
\begin{align}\label{FDPtW}
\widehat{\rm FDP}_{\tilde{W}}(t)\approx {\rm FDP}_{\tilde{W}}(t)\left(1+\frac{5t^3\bar{\kappa}}{18n}\right)+O_p\left(\sqrt{\xi_{n,p}/n}+\beta_p^{-1/2}\right).
\end{align}
The difference in the asymptotic expansions of $\widehat{\rm FDP}_W(t)$ and $\widehat{\rm FDP}_{\tilde{W}}(t)$ is due to the different large-deviation probabilities of the $W_j$ and $\tW_j$. 
This difference immediately suggests a ``refined" threshold as,
\begin{align}\label{lre}
L_{\rm refined}=\inf\left\{t>0:\frac{\#\{j: W_j\leq -t\}}{\#\{j:W_j\geq t\}\vee 1}\left\{1-\frac{4}{9}\theta(t)\right\}\leq \alpha\right\},
\end{align}
where
\[
\theta(t)=\frac{\left(\#\{j: W_j\leq -t\}-\#\{j: W_j\geq t\}\right)-\left(\#\{j: \tW_j\leq -t\}-\#\{j: \tW_j\geq t\}\right)}{\#\{j: W_j\leq -t\}\vee 1}.
\]
We show that $\theta(t)\approx -\frac{t^3\bar{\kappa}}{2n}$ under certain conditions, and consequently using $L_{\rm refined}$ is capable of eliminating the effect of the term $\frac{2t^3\bar{\kappa}}{9n}$ in (\ref{FDPW}).

The next theorem demonstrates that the refined procedure has better convergence rate in certain circumstance.
Basically, we restrict our attention to the sparse case, say $p_1/p_0\to 0$, such like $p_1=p^\eta$ for $0<\eta< 1$.
This is because the term of order $t^3/n$ only matters when $t$ is large. From the proof of Theorem \ref{thm1} we see that
$L\lesssim 2\log\{p/(\beta_p\alpha)\}$. In other words, only if $p_1$  or $\beta_p$ is small, the tail approximation of
$\widehat{\rm FDP}_{{W}}(t)$ to ${\rm FDP}_{{W}}(t)$ would be important.

\begin{thm}\label{thm2}
Assume $X_1,\ldots,X_p$ are independent each other. Suppose Assumptions \ref{moment}-(ii), \ref{signall}, $p=o\{\exp(n^{1/3})\}$, $p_1=p^{\eta}$ and $(p_1/\beta_p-1)\log^2(p/\beta_p)=O(1)$ hold. For any $\alpha\in(0,1)$, the FDP of the refined RESS method satisfies ${\rm FDP}_{\rm refined}\leq \alpha+O_p\left\{\sqrt{\xi_{n,p}/n}+n^{-1}(\log p)^{2}+\beta_p^{-1/2}\right\}$.
\end{thm}
In this theorem, the condition $(p_1/\beta_p-1)\log^2(p/\beta_p)=O(1)$ implies that the number of the signals we can identify dominates the number of those very weak signals. The RESS method with the refined threshold $L_{\rm refined}$ has the same convergence rate as the bootstrap calibration.
Simultaneous testing of many hypotheses allows us to construct a ``data-driven" correction of skewness without resampling. Thus, in some sense, large-scale $t$ tests with ultrahigh-dimension may not be considered as a ``curse'' but a ``blessing'' in our problem.

We summarize the refined RESS procedure as follows.

{\it Reflection via Sampling Splitting Procedure (RESS)}
\begin{itemize}
  \item\textit{Step 1:} Randomly split the data into two parts with equal size. Compute $\bar{X}_{kj}$ and $s_{kj}^2$ for $k=1,2$ and $j=1,\ldots,p$;
  \item \textit{Step 2:} Obtain $W_{j}$ and $\tilde{W}_{j}$ for $j=1,\ldots,p$.  
 Compute $\theta(t)$ for $t\in\{|W_j|\}_{j=1}^p$;
 \item \textit{Step 3:} Find the threshold $L_{\rm refined}$ by (\ref{lre}) and reject $\bH_{0j}$ if $W_j\geq L_{\rm refined}$.
\end{itemize}
The total computation complexity is of order $O(np+p\log p)$ and the procedure can be easily implemented even without high-level programming language. The R and Excel codes are available upon request.

We want to make some remarks on the use of $\tilde{W}_j$. As can be seen from (\ref{FDPtW}),  $\widehat{\FDP}_{\tilde{W}}(t)$ is an overestimate of the true FDP, and therefore yields a slightly more conservative procedure. In practice, if the computation is our major concern, using the RESS procedure with $\tilde{W}_j$ could be a safe choice.


\subsubsection{Dependent case}

We establish theoretical properties of the dependent $X$ case. The first result
is a direct extension of Proposition \ref{pro2}.
\begin{pro}\label{pro3}
Assume that $n_t\geq 4$. For any $\alpha\in(0,1)$, the RESS method satisfies
\[
\mbox{\rm FDR}\leq \min_{\epsilon\geq 0}\left\{\alpha(1+4\epsilon)+\Pr\left(\max_{j\in\I_0}\Delta_j'>\epsilon\right)\right\},
\]
where $\Delta_j'=\left|\Pr(W_j>0\mid |W_j|,{\bf W}_{-j})-1/2\right|$ and ${\bf W}_{-j}=(W_1,\ldots,W_p)^{\top}\setminus W_j$.
\end{pro}
Again, $\Delta_j'$ quantifies the effect of both the asymmetry of $X_j$ and the dependence between $W_j$ and ${\bf W}_{-j}$ on the FDR.

To achieve asymptotical FDR control, the following condition on the dependence structure is imposed.
\begin{assu}\label{corr}
(Correlation) For each $X_j$, assume that the number of variables $X_k$ that
are dependent with $X_j$ is no more than $r_p=o(\beta_p)$.
\end{assu}

This assumption implies that $X_j$ is independent with the other $p-r_p$ variables.
This is certainly not the weakest condition, but is adopted here to simplify the proof.

Let $\zeta_p=(r_p/\beta_p)^{1/2}$. The following theorem is parallel with Theorems \ref{thm1}-\ref{thm2}.
\begin{thm}\label{thm3} Suppose Assumptions \ref{signall}, \ref{corr} and $p=o\{\exp(n^{1/3})\}$ hold. \begin{itemize}
\item[(i)] If Assumption \ref{moment}-(i) hold, then for any $\alpha\in(0,1)$, the FDP of the RESS method satisfies
\begin{align*}
{\rm FDP}\leq\alpha+O_p\left\{\sqrt{\xi_{n,p}/n}+n^{-1}(\log p)^{3}\max_{j\in \I_0}\kappa_j^2+\zeta_p\right\},
\end{align*}
and $\mathop{\lim\sup}_{n\to\infty}{\rm FDR}\leq \alpha$.
\item[(ii)] If Assumption \ref{moment}-(ii), $p_1=o(p)$ and $(p_1/\beta_p-1)\log^2(p/\beta_p)=O(1)$ hold, then for any $\alpha\in(0,1)$, the FDP of the refined RESS method satisfies ${\rm FDP}_{\rm refined}\leq\alpha+O_p\left\{\sqrt{\xi_{n,p}/n}+n^{-1}(\log p)^{2}+\zeta_p\right\}$.
\end{itemize}
\end{thm}

This theorem implies that the RESS method remains valid asymptotically for weak dependence.
Comparing Theorem \ref{thm3} with Theorems \ref{thm1}-\ref{thm2}, the main difference lies on the  convergence rates of $\beta_p^{-1/2}$ and $\zeta_p$; the latter one is asymptotically larger. This can be understood because the approximation of the empirical distribution to the population one is expected to have slower rate for dependent summation.
When $X_j$ is independent with all the other variables $X_k$ or only dependent with finite number of $X_k$, then $r_p=O(1)$. In this situation, Theorem \ref{thm3} reduces to Theorems \ref{thm1}-\ref{thm2}.

\remark 3 In our discussion, we restrict $p=o\{\exp(n^{1/3})\}$ which facilitates our technical derivation. In \cite{liu2014phase}, a faster rate of $p$ is allowed, $p=o\{\exp(n^{1/2})\}$. However, we also notice that the bootstrap method proposed by \cite{liu2014phase} jointly uses the empirical distribution of $T_j^*,j=1,\ldots p$, where $T_j^*$ is the bootstrap statistic for the $j$th variable. This implies that its computation complexity is of order $O(p^2B+pnB)$ and it also requires $O(pB)$ storage, where $B$ is the number of bootstrap replications. For the commonly used bootstrap, say to approximate the distribution of $T_j$ individually by resampling \citep{FHY2007}, $p=o\{\exp(n^{1/2})\}$ is achieved if the replication of bootstrap is of order $p^2$. Though our theoretical results only allow $p=o\{\exp(n^{1/3})\}$, we conjecture that similar results also hold when $p=o\{\exp(n^{1/2})\}$ if more stringent conditions were imposed. Encouragingly, our extensive simulation results show that the refined procedure could work at least as good as the bootstrap method in terms of FDR control, even when $p$ is super-large. \hfill$\Box$

\section{Extensions}

In this section, we discuss two generalizations of our RESS procedure.

\subsection{One-sided alternatives}

In certain situations, we are interested in the case with one-sided
alternatives, say without loss of generality, $\mu_j<0$ for all $j\in\I_1$. We
may modify the RESS by ruling out the cells with $T_{1j}>0$ and $T_{2j}>0$
from the set $\{j: W_j\geq L\}$ to improve the power. To be more specific, the
threshold in (\ref{th}) is modified by
\[
\tilde{L}=\inf\left\{t>0:\frac{\#\{j: W_j\leq -t\}-\#\{j: W_j\geq t, T_{1j}>0,T_{2j}>0\}}{\#\{j:W_j\geq t,T_{1j}<0,T_{2j}<0\}\vee 1}\leq \alpha\right\},
\]
and we reject $\bH_{0j}:\mu_j\geq 0$ when $W_j\geq \tilde{L},T_{1j}<0,T_{2j}<0$. We have the following result.
\begin{coro}\label{coro1}
Consider the one-sided hypotheses that $\mu_j<0$ for all $j\in\I_1$.  If the conditions in Theorem \ref{thm3}-(i) hold,  then for any $\alpha\in(0,1)$, the FDP of the  RESS method with the threshold $\tilde{L}$ satisfies \begin{align*}
{\rm FDP}(\tilde{L})\leq\alpha+O_p\left\{\sqrt{\xi_{n,p}/n}+n^{-1}(\log p)^{3}\max_{j\in \I_0}\kappa_j^2+\zeta_p^{-1}\right\},
\end{align*}
and $\mathop{\lim\sup}_{n\to\infty}{\rm FDR}(\tilde{L})\leq \alpha$.
\end{coro}

By using the results in \cite{delaigle2011robustness}, the convergence rate of
the normal calibration is
\begin{align*}
 {\rm FDP}_{\Phi}&\lesssim\alpha+O_p\left\{n^{-1/2}(\log p)^{3/2}\max_{j\in \I_0}\kappa_j^2\right\}.
\end{align*}
Comparing this with Corollary \ref{coro1}, we see that the RESS strategy has
removed the skewness term that describes first-order inaccuracies of the
standard normal approximation. This important property is due to the fact that
$W_j$ is more symmetric than the $t$ statistic. See the proof of Theorem
\ref{thm1} for details. The refined RESS procedure, which also enjoys the
second-order accuracy, can be defined similarly to $L_{\rm refined}$ but we do
not discuss it in detail.

\subsection{Two-sample problem}

We next extend the RESS procedure to two-sample problem. Assume there is
another random sample $\mathcal{G}=(Z_{i1}, \cdots, Z_{ip})_{i=1}^{2m}$
distributed from $(Z_1,\cdots, Z_p)$. The population mean vectors of
$(X_1,\cdots, X_p)$ and $(Z_1,\cdots, Z_p)$ are $(\mu^X_{1},\cdots,\mu^X_p)$
and $(\mu^Z_{1},\cdots,\mu^Z_p)$, respectively. We aim to carry out $p$
two-sample tests, that is, $\bH_{0j}: \mu^{X}_j=\mu^{Z}_j$ versus $\bH_{1j}:
\mu^{X}_j\neq \mu^{Z}_j$, for $j=1, \ldots, p$. The RESS procedure can be
readily generalized to this two-sample problem as follows.

Firstly, similar to the sample splitting of $\D_{1}$ and $\D_{2}$, the data $\mathcal{G}$ are also splitted randomly into two disjoint groups $\mathcal{G}_{1}$ and $\mathcal{G}_{2}$ with equal size $m$. Based on $\D_k$ and $\mathcal{G}_k$, two-sample $t$-test $T_{kj}$ is defined as follows: 
$$T_{kj}=\frac{\bar{X}_{kj}-\bar{Z}_{kj}}{\sqrt{s^2_{X_{kj}}/n+s^2_{Z_{kj}}/m}}~~\text{for }j=1,\ldots,p~~{\text{and }}k=1,2.$$
Here $\bar{X}_{kj}$ and $\bar{Z}_{kj}$ are the sample means of $X_j$ and $Z_j$, while $s^2_{X_{kj}}$ and $s^2_{Z_{kj}}$ are the sample variances of $X_j$ and $Z_j$.
Finally, define $W_j=T_{1j}T_{2j}$. The threshold $L$ can be then selected similarly as in (\ref{th}), and $\bH_{0j}$ will be rejected when $W_j\geq L$.

To establish the FDR control result, the Assumptions 2 and 3 are modified as follows.
\begin{assu}\label{signallnew}
(Signals) For a large $C$, $\beta_p\equiv\ca\left\{j: |\mu^X_j-\mu^Z_j|\geq C\sqrt{\log p/n}\right\}\to \infty.$
\end{assu}

\begin{assu}\label{corrnew}
(Correlation) For each $X_j$, assume that the number of variables $X_k$ that
are dependent with $X_j$ is no more than $r_p=o(\beta_p)$. The same
assumption is imposed on $(Z_1,\cdots, Z_p)$.
\end{assu}

By Theorem 2.4 in \cite{chang2016cramer} and the proofs for Theorem 3, we have the following result.

\begin{coro}
Suppose Assumption 1-(i), Assumptions \ref{signallnew} and \ref{corrnew} and $p=o(\exp(n^{1/3}))$ hold,  then for any $\alpha\in(0,1)$, the FDP of the  RESS method for the two-sample problem satisfies \begin{align*}
{\rm FDP}(L)\leq\alpha+O_p\left\{\sqrt{\xi_{n,p}/n}+n^{-1}(\log p)^{3}+\zeta_p^{-1}\right\},
\end{align*}
and $\mathop{\lim\sup}_{n\to\infty}{\rm FDR}\leq \alpha$.
\end{coro}

\section{Numerical results}

\subsection{Simulation comparison}

We evaluate the performance of our proposed RESS procedure on simulated data sets and compare the FDR and true positive rate (TPR) with other existing techniques. All the results are obtained from 200 replication simulations.

\subsubsection{Model and benchmarks}

We set the model as $X_{ij}=\mu_j+\varepsilon_{ij},~~1\leq i\leq n_t,~~1\leq
j\leq p$, where the alternative signal using $\mu_j=\delta_j\sqrt{\log p/n_t}$ with
$\delta_j\sim \text{Unif}(-1.5,-1)$ or $\delta_j\sim \text{Unif}(1,1.5)$ for
$j\in\mathcal{I}_1$. The random error is generated by the autoregression model
$\varepsilon_{ij}=\rho\varepsilon_{i,j-1}+\epsilon_{ij}$, where $\rho\in[0,1)$
and $\epsilon_{ij}$'s are i.i.d from three centered distributions: Student's $t$
with five degrees of freedom $t(5)$, exponential with rate one
($\text{exp}(1)-1$) and a mixed distribution which consists of
$\epsilon_{ij}\sim N(0,1)$ for $1\leq j\leq [p/3]$, $\epsilon_{ij}\sim t(5)$
for $[p/3]+1\leq j\leq [2p/3]$ and $\epsilon_{ij}\sim\text{exp}(1)-1$ for
$[2p/3]+1\leq j\leq p$. When $\rho=0$, the random errors are independent. We
consider the number of alternatives as $p_1=[\vartheta p]$ with
$\vartheta\in[0.01,0.15]$. 

The following three benchmarks are considered for comparison. The first one is the BH
procedure with the $p$-values obtained from the standard normal approximation (termed as
BH for simplicity). The other two are bootstrap-based approaches. Assume
$\{{\bm X}^*_{k1},\ldots {\bm X}^*_{kn_t}\}$, $1\leq k\leq B$ denotes bootstrap resamples
drawn by sampling randomly and $T^{*}_{kj}$ are Student's {\it{t}} test statistics
constructed from $\{X^{*}_{1kj}-\bar{X}_j,\ldots,X^{*}_{n_tkj}-\bar{X}_j\}$. One bootstrap
method is to estimate the $p$-values according to the bootstrap distribution individually
by $p_{j,BI}=1/B\sum_{k=1}^B\bI\{|T^{*}_{kj}|\geq|T_j|\}$ \citep{FHY2007} and another one
is to estimate the $p$-values with the average $p$ bootstrap distribution together by
$p_{j,BA}=1/{Bp}\sum_{k=1}^B\sum_{i=1}^p\bI\{|T^{*}_{ki}|\geq |T_j|\}$
\citep{liu2014phase}. We call these two as ``I-bootstrap"  and ``A-bootstrap",
respectively.  The A-bootstrap jointly estimates the distribution of $T_j$'s and thus can
be expected to have better performance. However, the computational complexity of
I-bootstrap is much lower; it is of order $O(npB)$ while that of A-bootstrap increases in
a quadratic rate of $p$. We take the number of bootstrap samples as $B=200$ in this
section.


\subsubsection{Results}

We compare the performance of our proposed RESS method
(termed as ``$\text{RESS}_0$") and the refined RESS method (``RESS")
in a range of settings, with the BH procedure, I-bootstrap and A-bootstrap,
and examine the effects of skewness, signal magnitude and correlation between variables.

\begin{table}
\caption{\small Comparison results of FDR, TPR and the average computing time under $p=5000, \rho=0$ and $\vartheta=0.05$.
Values in parentheses are standard deviations.} \label{tab1}\vspace{0.3cm}
\linespread{1.2}
\centering
{\small\renewcommand{\arraystretch}{1}
\begin{tabular}{cccccccccccccccccc}
\hline
Errors&&  \multicolumn{3}{c}{$n_t=50$} && \multicolumn{3}{c}{$n_t=100$} \\
& Method & \tc{FDR(\%)} & \tc{TPR(\%)} & \tc{Time(s)} && \tc{FDR(\%)} & \tc{TPR(\%)} & \tc{Time(s)} \\
\hline
    &   RESS  &   $17.8_{\tiny{(4.6)}}$    &   $51.6_{\tiny{(5.5)}}$    &   1.75  &&  $18.5_{\tiny{(5.0)}}$    &   $51.6_{\tiny{(5.8)}}$    & 1.84     \\
    &   $\text{RESS}_0$   &   $19.7_{\tiny{(4.4)}}$    &   $53.8_{\tiny{(5.4)}}$    &   1.66  &&  $19.9_{\tiny{(4.9)}}$    &   $53.1_{\tiny{(5.6)}}$    & 1.74    \\
$t(5)$  &   BH  &   $16.6_{\tiny{(3.9)}}$  &   $52.1_{\tiny{(4.4)}}$  & 1.60    &&  $17.7_{\tiny{(3.2)}}$  &   $52.9_{\tiny{(4.5)}}$  & 1.60 \\
    &   A-bootstrap &   $12.8_{\tiny{(3.4)}}$  &   $47.6_{\tiny{(4.8)}}$  & 404 &&  $16.3_{\tiny{(3.1)}}$  &   $51.4_{\tiny{(4.5)}}$  & 400 \\
    &   I-bootstrap &   $22.3_{\tiny{(2.6)}}$  &   $53.2_{\tiny{(3.3)}}$  & 314 &&  $23.6_{\tiny{(3.0)}}$  &   $54.9_{\tiny{(3.3)}}$  & 311 \\
\\
    &   RESS  &   $19.4_{\tiny{(5.0)}}$    &   $63.1_{\tiny{(10.1)}}$    &   1.76 &&  $20.8_{\tiny{(4.9)}}$    &   $75.2_{\tiny{(5.8)}}$    & 1.85     \\
    &   $\text{RESS}_0$   &   $32.2_{\tiny{(3.8)}}$    &   $81.5_{\tiny{(3.6)}}$    &   1.66  &&  $27.3_{\tiny{(3.9)}}$    &   $81.4_{\tiny{(2.8)}}$    & 1.76 \\
exp(1)  &   BH  &   $37.2_{\tiny{(2.7)}}$  &   $85.5_{\tiny{(2.9)}}$  & 1.70    &&  $30.2_{\tiny{(2.7)}}$ &   $84.4_{\tiny{(2.5)}}$  &  1.63    \\
    &   A-bootstrap &   $18.8_{\tiny{(3.5)}}$  &   $62.8_{\tiny{(5.1)}}$  & 424 &&  $19.7_{\tiny{(2.8)}}$  &   $76.0_{\tiny{(3.4)}}$  & 403 \\
    &   I-bootstrap &   $32.2_{\tiny{(3.5)}}$  &   $70.2_{\tiny{(4.6)}}$  & 332 &&  $28.2_{\tiny{(2.6)}}$  &   $77.9_{\tiny{(2.6)}}$  & 317 \\
    \\
    &   RESS  &   $20.7_{\tiny{(3.9)}}$    &   $67.7_{\tiny{(5.3)}}$    &   1.85  &&  $20.6_{\tiny{(4.0)}}$    &   $70.0_{\tiny{(4.4)}}$    & 2.17     \\
    &   $\text{RESS}_0$   &   $25.1_{\tiny{(3.9)}}$    &   $72.2_{\tiny{(4.0)}}$    &   1.75  &&  $23.1_{\tiny{(3.7)}}$    &   $72.1_{\tiny{(3.8)}}$    & 2.06    \\
Mixed   &   BH  &   $26.2_{\tiny{(2.9)}}$  &   $74.0_{\tiny{(3.1)}}$  & 1.65    &&  $23.7_{\tiny{(3.1)}}$  &   $74.3_{\tiny{(3.2)}}$  & 1.64    \\
    &   A-bootstrap &   $17.5_{\tiny{(3.0)}}$  &   $65.1_{\tiny{(4.4)}}$  & 405 &&  $19.0_{\tiny{(2.9)}}$  &   $70.2_{\tiny{(3.5)}}$  & 409 \\
    &   I-bootstrap &   $25.7_{\tiny{(4.4)}}$  &   $67.3_{\tiny{(5.6)}}$  & 314 &&  $26.9_{\tiny{(4.1)}}$  &   $72.8_{\tiny{(4.6)}}$  & 319 \\
\hline
\end{tabular}}
\end{table}

Firstly, we set $p=5000$ and $\vartheta=0.05$. The full sample size $n_t$ is taken to be
$50$ or $100$, and the target FDR level $\alpha$ is set as $0.2$.
Table 1 displays the estimated FDR, TPR and average
computation time obtained by each method under three different
error settings with $\rho=0$.  For the symmetric distribution $t(5)$,
 the $\text{RESS}_0$ is able to deliver a quite accurate control
but that is not the case for the other skewed errors. It performs better than the BH in
terms of FDR control, but has a slightly disadvantage over the I-bootstrap. This is
consistent with our theoretical analysis in Section 2.2.1.
{\color{black}Actually, from Theorem 1 and Theorem 3-(i), when $p$ is very large, the skewness still has non-ignorable effect on the FDR control of the $\text{RESS}_0$ method.} In contrast, we observe that
the FDR levels of $\text{RESS}$ are close to the nominal level under all the scenarios;
it is clearly more effective than the I-bootstrap, $\text{RESS}_0$ and BH, and the
difference is quite remarkable in some cases. Certainly, this is not surprising as it is
a data-driven method which has second-order accuracy justified in Section 2.2.2. The
power of RESS is also quite high compared to the other methods, revealing that its
detection ability is not largely compromised by data-splitting.

\begin{figure}[ht]\centering
\includegraphics[width=1\textwidth]{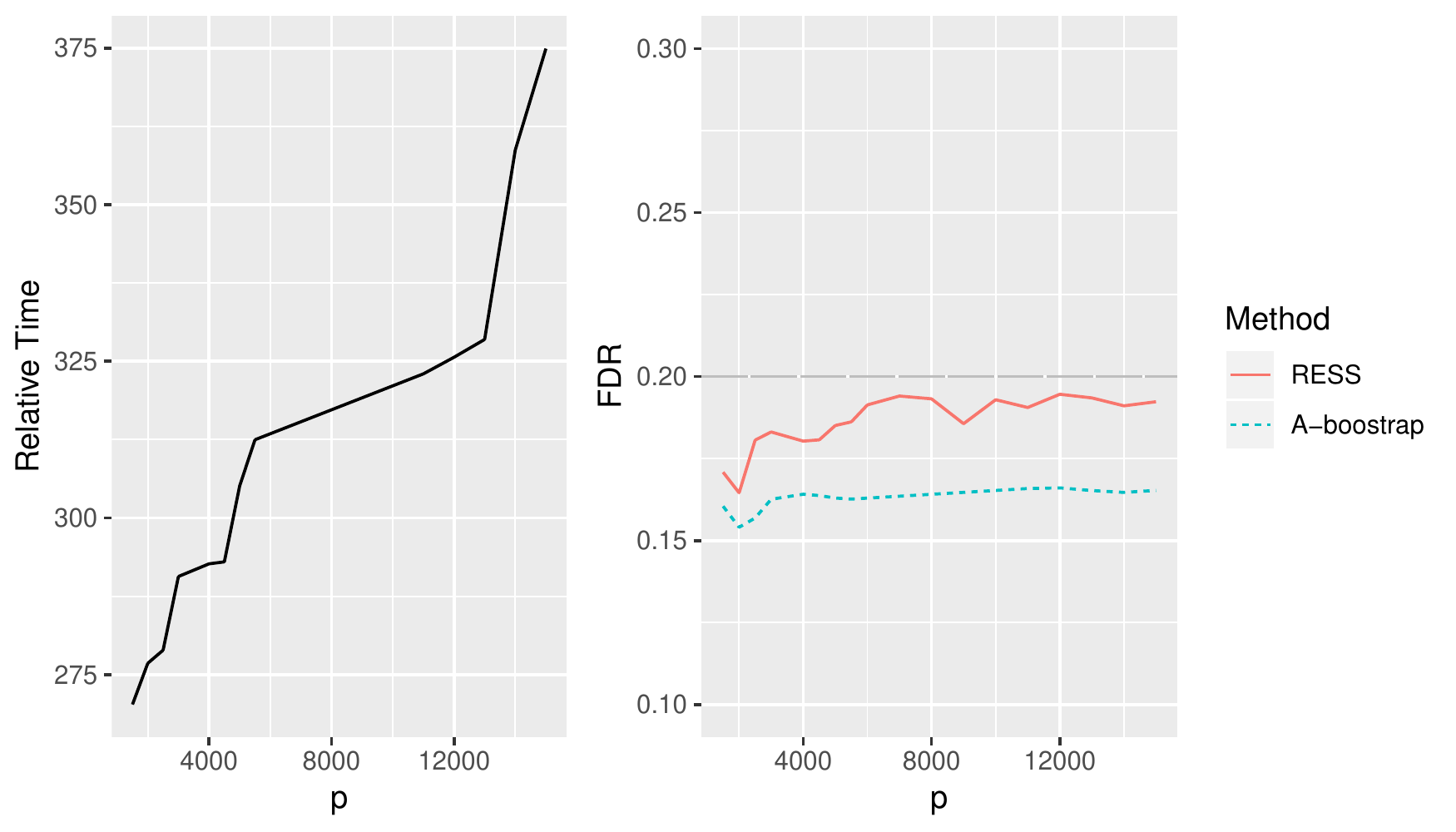}\vspace{-0.2cm}
\caption {\it\small Relative times and FDR curves against the number of tests $p$ when errors are distributed from $\text{t}(5)$ under $n_t=100$, $\rho=0$ and $\vartheta=0.05$; The grey dashed lines indicate the target FDR level.}\label{Fig:p}\vspace{-0.2cm}
\end{figure}

The I-bootstrap improves slightly the accuracy of FDR control under the $\text{exp}(1)$
and mixed cases over the normal calibration, but its FDRs are still considerably larger
than the nominal level. Note that the I-bootstrap calibration may need an extremely large
replications $B$, i.e. $p^2$, to achieve FDR control \citep{liu2014phase}, and therefore
does not perform well with commonly used $B=200$. The A-bootstrap method offers
comparable performances to our RESS method, though it tends to be a little conservative
under the $t(5)$ cases. Also, the A-bootstrap generally has smaller variations than
RESS due to the use of bootstrap for estimating an overall empirical null distribution.
Certainly, this gain comes partly from its computation-intensive feature; in most cases,
it requires more than 200 times computational time than the RESS. In fact, as mentioned
earlier, the computational complexities of A-bootstrap and RESS are of order $p^2B+pn_tB$
and $p\log p+pn_t$, respectively, and hence their relative computational time could
increase fast as $p$ increases. Figure \ref{Fig:p} depicts the relative time of
A-bootstrap to RESS and the FDR vaules under the $t(5)$ cases.

\begin{figure}[ht]\centering
\includegraphics[width=1\textwidth]{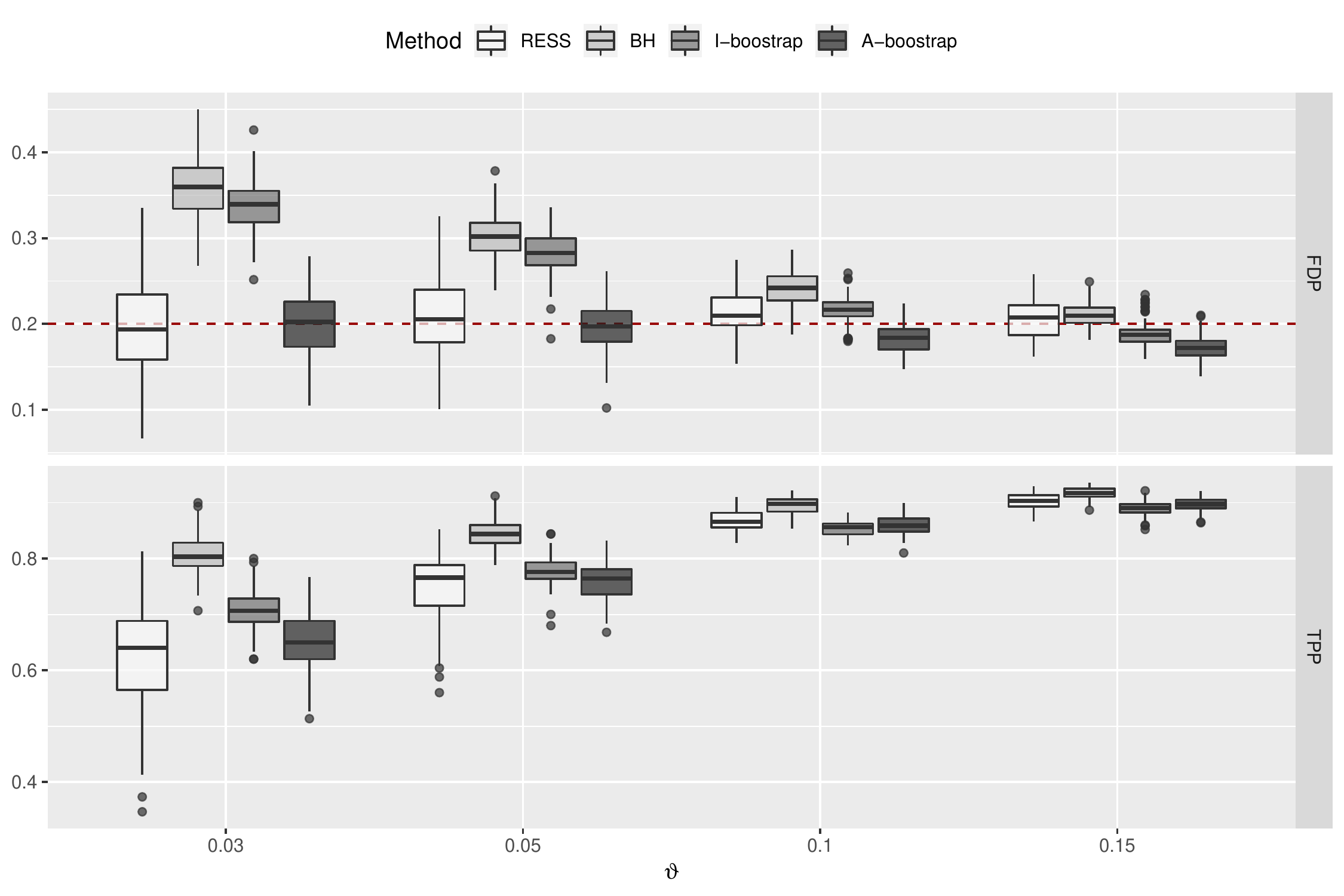}\vspace{-0.2cm}
\caption {\it\small Box-plots of FDP and true positive proportion (TPP) when errors are distributed from $\exp(1)$ under $n_t=100$, $p=5000$, and $\rho=0$; The red dashed lines indicate the target FDR level.}\label{Fig:ratio}\vspace{-0.2cm}
\end{figure}

\begin{figure}[ht]\centering
\includegraphics[width=1\textwidth]{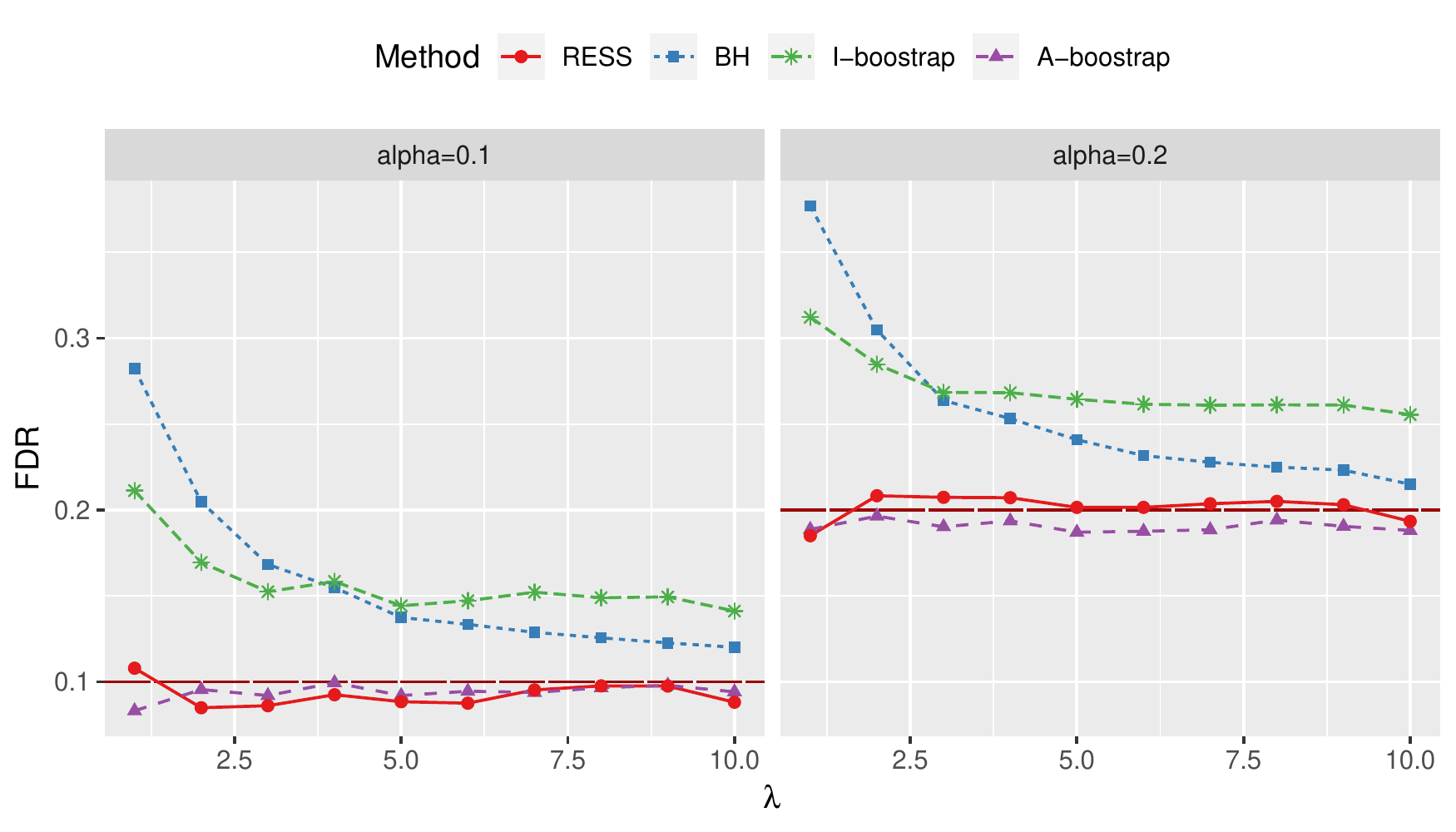}\vspace{-0.2cm}
\caption {\it\small FDR curves against the values of {$\lambda$ when errors are distributed from Gamma distribution $\Gamma(\lambda/2,2)$}
under $n_t=100$, $p=5000$, $\rho=0$ and $\vartheta=0.05$; The red dashed lines indicate the target FDR level.}\label{Fig:chi}\vspace{-0.2cm}
\end{figure}

Next, we examine the effect of the skewness and the proportion of alternatives.  As we
have shown that the refined RESS performs usually better than the  ${\rm RESS}_0$, in
what follows we focus only on the refined RESS. With respect to skewness, we evaluate the performances of various methods by varying shape parameter $\lambda$ of Gamma distribution $\Gamma(\lambda/2,2)$, 
Figure \ref{Fig:chi} depicts the estimated FDR curves against the
values of $\lambda$ when $p=5000$, $n_t=100$ and $\vartheta=0.05$ for $\alpha=0.1$ and
$\alpha=0.2$. It can be seen that the refined RESS successfully controls FDR in an acceptable
range of the target level no matter the magnitude of the skewness. Figure \ref{Fig:ratio}
shows the boxplots of FDP and powers under $\vartheta=0.03,0.05,0.1$ and $0.15$ when the
errors are generated from $\text{exp}(1)$. The refined RESS has stable FDP close to
the nominal level $\alpha=0.2$, even when $\vartheta$ is small. Similar results with
other error distributions and signal magnitudes can be found in the Supplementary
Material.

\begin{figure}[ht]\centering
\includegraphics[width=1\textwidth]{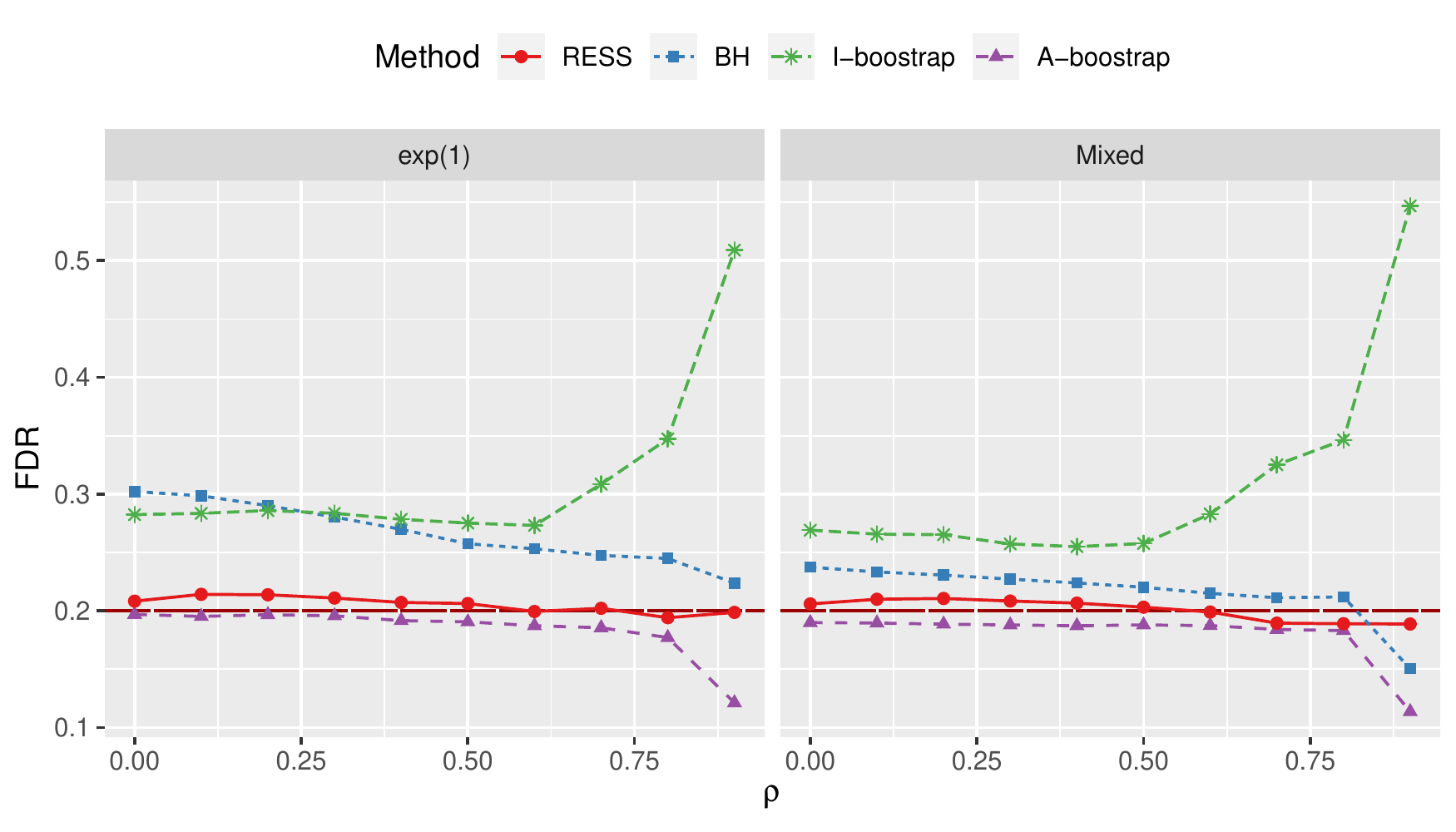}\vspace{-0.2cm}
\caption {\it \small FDR curves against the correlation when
errors are distributed from $\text{exp}(1)$ and mixed one under $n_t=100$,
$p=5000$ and $\vartheta=0.05$; The red dashed lines indicate the target FDR level.}
\label{Fig:rho}\vspace{-0.2cm}
\end{figure}

Finally, we turn to investigate the effect of the correlation level $\rho$. Figure
\ref{Fig:rho} shows the FDR curves of four methods against the values of $\rho$ when the
errors are generated from the $\text{exp}(1)$ and mixed distributions. Again, it can be seen that
the refined RESS and A-bootstrap result in a reasonably good FDR control in most
situations, even when $\rho$ is as large as 0.8. This concurs with our asymptotic
justification that the refined RESS method is still effective provided that ${X}$
satisfies certain weak dependence structure.




\subsection{A real-data example}

We next illustrate the proposed refined RESS procedure by an empirical analysis of the
acute lymphoblastic leukemia (ALL) data, which consists of 12,256 gene probe sets for 128
adult patients enrolled in the Italian GIMEMA multi center clinical trial 0496. It is
known that malignant cells in B-lineage ALL often have genetic abnormalities, which have
a significant impact on the clinical course of the disease. Specifically, the molecular
heterogeneity of the B-lineage ALL is well established as BCR/ABL, ALL1/AF4, E2A/PBX1 and
NEG and the gene expression profiles of groups BCR/ABL and NEG are more similar to each
other than to the others. In our analysis, we consider a sub-dataset of 79 B-lineage
units with 37 BCR/ABL mutation and 42 NEG and use the traditional two-sample $t$-test to
examine which probe sets are differentially expressed. The dataset was previously studied
by \citet{chiaretti2005gene} and \citet{Bourgon_etal_2010} and is available at {\it
http://www.bioconductor.org}.

\begin{figure}[ht]\centering
\includegraphics[width=1\textwidth]{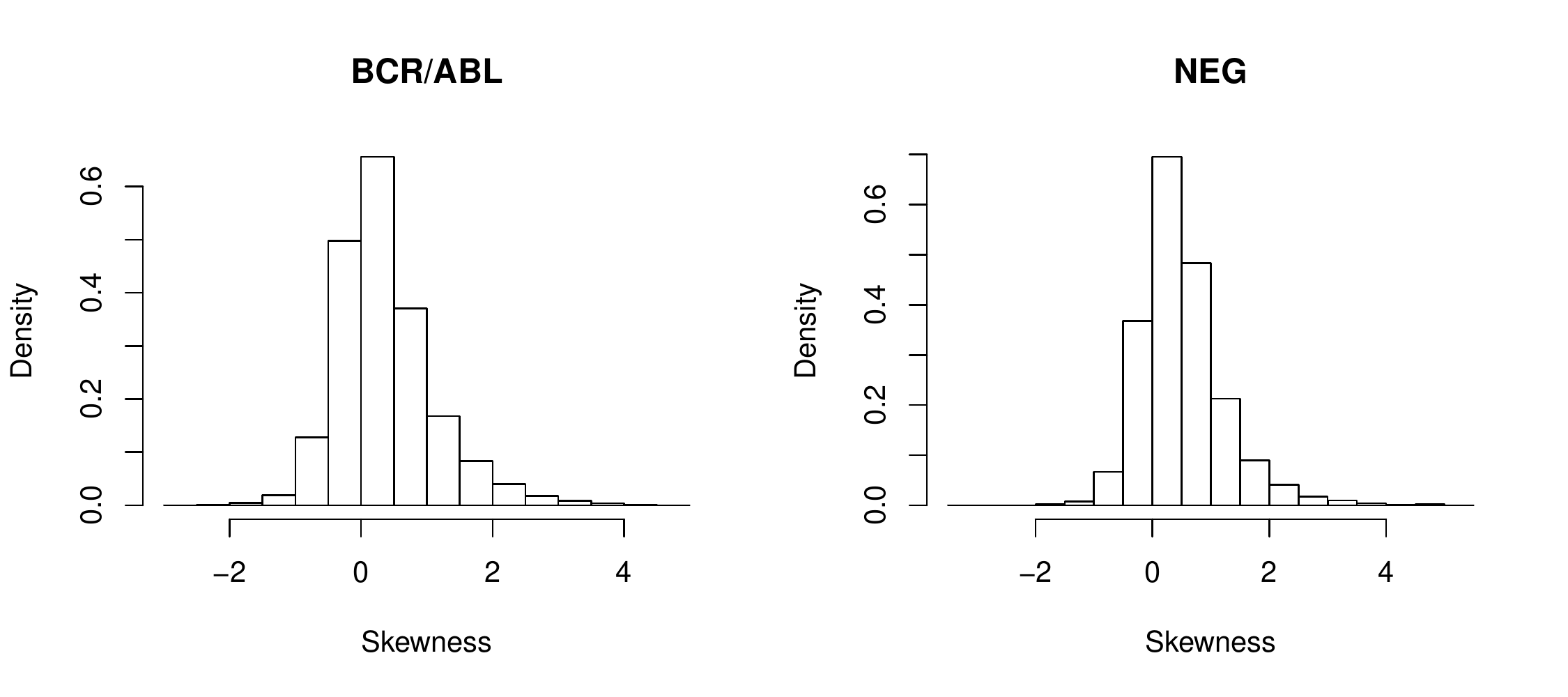}\vspace{-0.2cm}
\caption {\it \small Density histogram of the skewness of the $p=12,256$
genes for BCR/ABL and NEG, respectively.}
\label{Fig:skeALL}\vspace{-0.2cm}
\end{figure}

\begin{figure}[ht]\centering
\includegraphics[width=0.8\textwidth]{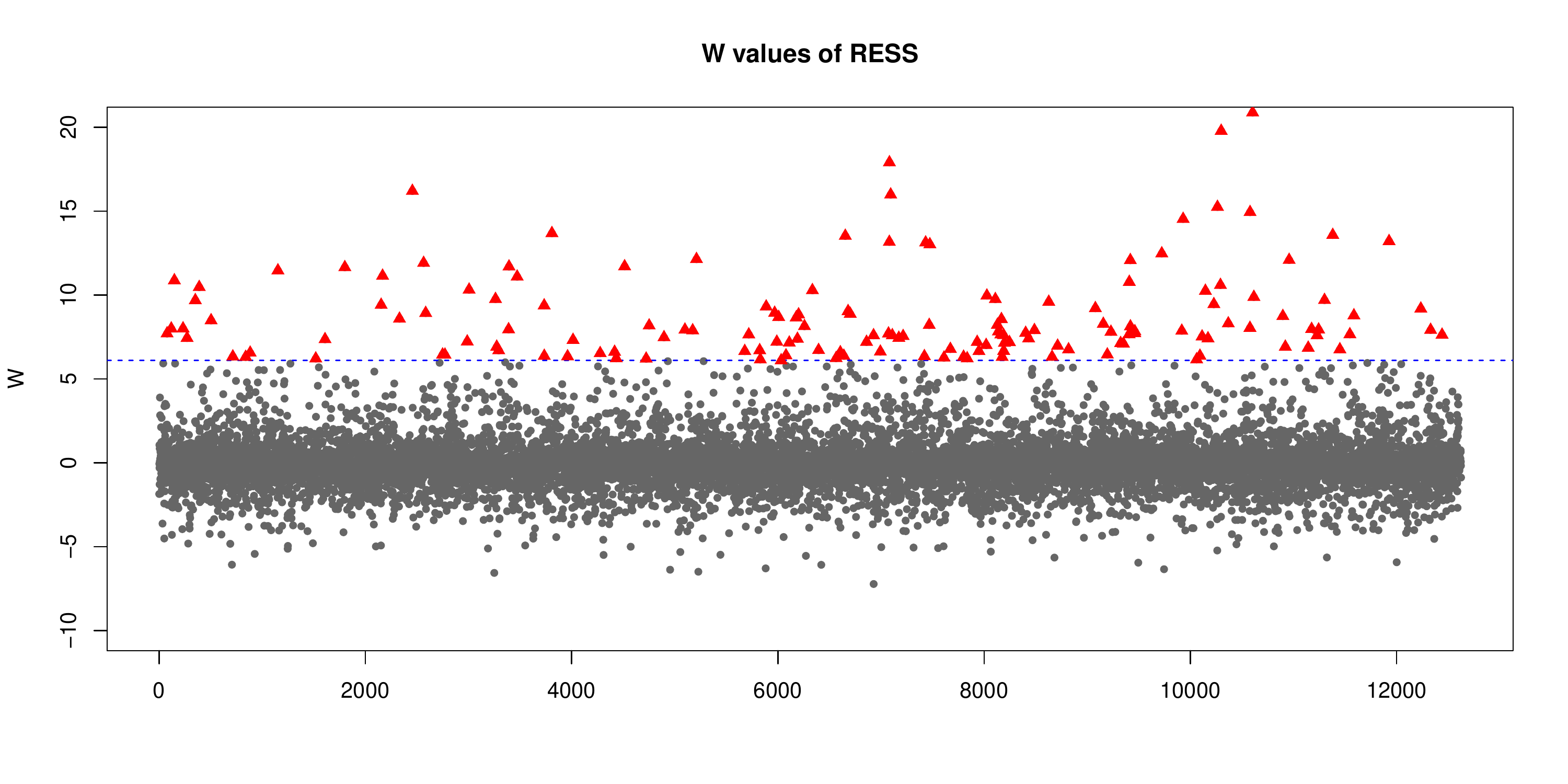}\vspace{-0.2cm}
\caption {\it \small Scatter plot of the $W_j$ for our raw RESS with the
red triangles and blue dashed line denoting selected differentially expressed probe sets and
the threshold under $\alpha=0.05$.}\label{Fig:WvalALL}\vspace{-0.2cm}
\end{figure}

Here, we consider the extension of our proposed refined RESS in two-sample case and
compare it with BH procedure, A-bootstrap and I-bootrstrap over a wide range of
significant levels. Both of the bootstrap sampling are repeated 200 times.
Table 2 summarizes the number of probe sets differentially
expressed between BCR/ABL and NEG. The BH procedure with the normal calibration
tends to reject surprising more genes than RESS and A-bootstrap for various
significance levels. In fact, the normality approximation seems to be violated
for many of the genes as some skewness values largely deviate from zero in
Figure \ref{Fig:skeALL}. As noted earlier, the I-bootstrap needs an extremely
large replications, i.e. $p^2$ to improve the accuracy and thus leads to the
large number of rejections with $B=200$. Figure \ref{Fig:WvalALL} presents the
scatterplot of RESS's statistics $W_j$. We can observe that all selected
probe sets (red triangles) have large values of $W_j$, while the unselected ones
(black dots) are roughly symmetric across the horizontal lines. With the
data-driven threshold, the number of differentially expressed probe sets based
on our RESS appears more reasonable. The results of A-bootstrap coincide with that of RESS as expected.

\begin{table}
\caption{\small Numbers of differentially expressed probe sets under various significant levels.} \label{tab2}\vspace{0.3cm}
\centering {\small
\renewcommand{\arraystretch}{1}\tabcolsep 9pt
\begin{tabular}{cccccccccccccccccc} \hline
$\alpha$ & \multicolumn{1}{c}{RESS} & \multicolumn{1}{c}{BH} & \multicolumn{1}{c}{I-bootstrap}   &   \multicolumn{1}{c}{A-bootstrap}\\
\hline
0.05 & 157  &   163 &   326 &   141\\
0.1 & 221   &   238 &   326 &   222\\
0.15    &   280 &   334 &   476 &   310 \\
0.2 &   333 &   414 &   476 &   397 \\
\hline
\end{tabular}}
\end{table}

\section{Concluding Remarks}

In this paper, we have proposed a multiple testing procedure, RESS,
that controls FDR in the large-scale $t$ setting and offers high power to discover
true signals. We give theoretical results showing
that the proposed method maintains FDR control under mild conditions.
The empirical performance of the refined RESS demonstrates excellent FDR
control and reasonable power in comparison to other methods such as the bootstrap
or normal calibrations. The ideas of the RESS procedure may be extended for controlling other
rates such as per family error rate.

\appendix

\setcounter{equation}{0}
\renewcommand{\theequation} {A.\arabic{equation}}
\def\thelemma{A.\arabic{lemma}}
\def\thepro{A.\arabic{pro}}

\subsection*{Appendix: Proofs}

\subsubsection*{Appendix A: Lemmas}

Before we present the proofs of the theorems, we first state several lemmas whose proofs
can be found in the Supplementary Material. A few well-known theorems to be repeatedly
used are also presented in the Supplementary Material. The first lemma characterizes the
closeness between $\Pr(W_j\geq t)$ and ${\Pr(W_j\leq -t)}$, which plays an important role
in the proof. For simplicity, we suppress the dependence of $T_{kj}$ on $j$ which should
not cause any confusion. Let $\tilde{\Phi}(x)=1-\Phi(x)$, $G(t)=p_0^{-1}\sum_{j\in
\I_0}\Pr(W_j\geq t)$, $G_{-}(t)=p_0^{-1}\sum_{j\in \I_0}\Pr(W_j\leq-t)$ and
$G^{-1}(y)=\inf\{t\geq 0: G(t)\leq y\}$ for $0\leq y\leq 1$.
\begin{lemma}\label{klem1} 
Suppose Assumption \ref{moment}-(i) hold. For any $0\leq t\leq C\log p$ with $C>0$,
\[
\frac{\Pr(T_1T_2\geq t)}{\Pr(T_1T_2\leq-t)}-1=O(a_{nt}n^{-1/2})+\frac{2t^3\kappa^2}{9n},
\]
where $a_{nt}=({2t}+{\log n})^{1/2}$.
\end{lemma}

The second lemma characterizes the closeness between $\Pr(\tilde W_j\geq t)$ and ${\Pr(\tilde W_j\leq -t)}$.
\begin{lemma}\label{klmwt}
Suppose Assumption \ref{moment}-(ii) hold. For sufficiently large $t$ satisfying $t\lesssim\log p$,
\[
\frac{\Pr(\tilde{W}\geq t)}{\Pr(\tilde{W}\leq -t)}-1=O(a_{nt}n^{-1/2}+t^{-1/2})-\frac{5t^3\kappa^2}{18n}.
\]
\end{lemma}

The next lemma establishes the uniform convergence of $\{{p_0G(t)}\}^{-1}{\sum_{j\in\I_0}\bI(W_j\geq t)}-1$.
\begin{lemma}\label{klm2}
Suppose Assumptions \ref{moment}-(i) and \ref{corr} hold. Then, for any $b_p\to\infty$ and $b_p=o(p)$,
\begin{align}
\sup_{0\leq t\leq G^{-1}(\alpha b_p/p)}\left|\frac{\sum_{j\in\I_0}\bI(W_j\geq t)}{p_0G(t)}-1\right|=O_p(\zeta_p),\label{mco1}\\
\sup_{0\leq t\leq G^{-1}_{-}(t)(\alpha b_p/p)}\left|\frac{\sum_{j\in\I_0}\bI(W_j\leq-t)}{p_0G_{-}(t)}-1\right|=O_p(\zeta_p)\label{mco2}.
\end{align}
\end{lemma}

The last one is the counterpart of Lemma \ref{klm2} for $\tilde{W}_j$.
\begin{lemma}\label{klem4}
Suppose Assumptions \ref{moment}-(ii) and \ref{corr} hold. Then, for any $b_p\to\infty$ and $b_p=o(p)$,
$(\ref{mco1})$ and (\ref{mco2}) holds if we replace ${W}_j$ and $G(t)$ with $\tilde{W}_j$ and $\tilde{G}(t)$, where $\tilde{G}(t)=p_0^{-1}\sum_{j\in \I_0}\Pr(\tilde{W}_j\geq t)$.
\end{lemma}

\subsubsection*{Appendix B: Proof of Theorems}

In fact, we show that under a weaker condition, i.e., Assumption \ref{corr}, the results in Theorems \ref{thm1}-\ref{thm2} hold similarly, and accordingly Theorem \ref{thm3} is also proved.

\noindent{\bf Proof of Theorem \ref{thm1}.} By definition, our test is equivalent to
reject $\bH_{0j}$ if $W_j\geq \wht$, where
\[
\wht=\inf\left\{t\geq 0: \sum_{j=1}^p\bI(W_j\leq-t)\leq \alpha\max\left(\sum_{j=1}^p\bI(W_j\geq t),1\right) \right\}.
\]

Let $\mathcal{A}$ be a subset of $\{1, 2,\ldots, p\}$ satisfying $\mathcal{A\subset}\left\{j: |\mu_j|/\sigma_j\geq 3\sqrt{\log p/n}\right\}$ and $b_p\equiv|\mathcal{A}|= \min(\sqrt{n\xi_{n,p}},\beta_p)$. By Assumption \ref{moment}-(i) and Markov's inequality, for any $\epsilon>0$, $\Pr(\max_{j \in\mathcal{A}}|s_{kj}^2/\sigma^2_j-1|\geq \epsilon)=O(\sqrt{\xi_{n,p}/n}),k=1,2$.
By Assumption \ref{signall}  and Lemma S.1, there exists some $c>2$,
\begin{align}
 &\Pr\left(\sum_{j=1}^p \bI(W_j\geq c\log p)\geq b_p/2\right)\nonumber\\
 &\geq \Pr\left(\sum_{j\in\mathcal{A}}\left\{\bI(T_{1j}T_{2j}\geq c\log p)-\Pr(T_{1j}T_{2j}\geq c\log p)\right\}\geq b_p/2-\sum_{j\in\mathcal{A}}\Pr(T_{1j}T_{2j}\geq c\log p)\right)\nonumber\\
 &\geq 1-\frac{\sum_{j\in\mathcal{A}}\Pr(T_{1j}T_{2j}\geq c\log p)}{\left\{\sum_{j\in\mathcal{A}}\Pr(T_{1j}T_{2j}\geq c\log p)-b_p/2\right\}^2}\nonumber\\
 &\geq 1-cb_p^{-1}\left\{1+o(1)\right\}=1+O(\sqrt{\xi_{n,p}/n}+b_p^{-1}),\label{tt}
\end{align}
where we use the fact that $\Pr_{j\in\mathcal{A}}(T_{1j}>\sqrt{c\log p})\geq \Pr_{j\in\I_0}(T_{1j}\geq-\sqrt{c\log p})\geq (1-p^{-1})\left\{1+o(1)\right\}$.

Define $\eta_p=\sqrt{2\log \{p/(b_p\alpha)\}}$.  We note that $G(t)\geq \alpha b_p/p$ implies that $t\leq \eta_p^2$ when $p$ is large. This is because
\begin{align*}
G(\eta^2_p)\leq 2\Pr(|T_1|>{\eta_p})
\leq 2\sqrt{2/\pi}\exp(-\eta^2_p/2)/\eta_p
\leq \alpha b_p/p\leq G(t),
\end{align*}
if $p$ is sufficiently large. This, together with Lemma \ref{klm2} and Eq. (\ref{tt}) implies that $L\leq \eta_p$ with probability at least $O(\sqrt{\xi_{n,p}/n}+b_p^{-1})$.

Therefore, by Lemma \ref{klem1}, we get
\begin{align}\label{fl1}
\frac{\sum_{j\in \I_0}\bI(W_j\geq \wht)}{\sum_{j\in \I_0}\bI(W_j\leq-\wht)}-1=O_p(\sqrt{\xi_{n,p}/n})+\frac{2\wht^3\max_{j\in\I_0}\kappa_j^2}{9n}\left\{1+O(a_{n\wht}/\sqrt{n}+\zeta_p)\right\}.
\end{align}
Write
\begin{align*}
{\rm FDP}&=\frac{\sum_{j\in\I_0}\bI\left(W_j\geq \wht\right)}{1\vee\sum_j\bI(W_j\geq \wht)}=\frac{\sum_{j}\bI\left(W_j\leq- \wht \right)}{1\vee\sum_j\bI(W_j\geq \wht)}\times\frac{\sum_{j\in\I_0}\bI\left(W_j\geq \wht\right)}{\sum_{j}\bI\left(W_j\leq-\wht\right)}\\
&\leq \alpha\times R(\wht).
\end{align*}
Note that $R(\wht)\leq {\sum_{j\in\I_0}\bI\left(W_j\geq \wht\right)}/{\sum_{j\in\I_0}\bI\left(W_j\leq -\wht\right)}$, and thus $\mathop{\lim\sup}_{n\to\infty} {\rm FDP}\leq\alpha$ in probability by (\ref{fl1}). Then, for any $\epsilon>0$,
\[
\FDR\leq (1+\epsilon)\alpha R(\wht)+\Pr\left(\FDP\geq (1+\epsilon)\alpha R(\wht)\right),
\]
from which the second part of this theorem is proved. \hfill$\Box$

\noindent{\bf Proof of Theorem \ref{thm2}.} We firstly establish a lower bound for $L$ so
that the condition in Lemma \ref{klmwt} is valid.  Notice
\begin{align*}
\sum_j\bI(W_j\leq-t)&\geq \sum_{j\in\I_0}\bI(W_j\leq-t) \approx p_0G_{-}(t),\\
\alpha\sum_j\bI(W_j\geq t)&\leq \alpha \left\{p_1+\sum_{j\in\I_0}\bI(W_j\geq t)\right\}\approx \alpha \left\{p_1+p_0G(t)\right\}.
\end{align*}
Hence, we can conclude that $L\gtrsim G^{-1}(\alpha p_1/\{p_0(1-\alpha)\})$.

For $\log\{p_0(1-\alpha)/(\alpha p_1)\}\lesssim t\lesssim\log \{p/(b_p\alpha)\}$,
\begin{align*}
\theta(t)=&\frac{\left\{\sum_{j}\bI(W_j\leq- t )-\sum_{j}\bI(\tW_j\leq- t )\right\}-\left\{\sum_{j}\bI(W_j\geq t )-\sum_{j}\bI(\tW_j\geq t )\right\}}{\sum_{j}\bI(W_j\leq -t)}\\
=&\frac{\sum_{j\in\I_0}\bI(W_j\leq -t)}{{1\vee\sum_j\bI(W_j\leq -t)}}\times\Bigg[\frac{\sum_{j\in\I_0}\left\{\bI(W_j\leq -t )-\bI(W_j\geq t )\right\}}{\sum_{j\in\I_0}\bI(W_j\leq -t)}+\frac{\sum_{j\in\I_0}\left\{\bI(\tW_j\geq t )-\bI(\tW_j\leq -t )\right\}}{\sum_{j\in\I_0}\bI(W_j\leq -t)}\\
&+\frac{\sum_{j\in\I_1}\left\{\bI(\tW_j\geq t )-\bI(W_j\geq t )\right\}}{\sum_{j\in\I_0}\bI(W_j\leq -t)}+\frac{\sum_{j\in\I_1}\left\{\bI(W_j\leq- t )-\bI(\tW_j\leq- t )\right\}}{\sum_{j\in\I_0}\bI(W_j\leq -t)}\Bigg]\\
\equiv& \frac{\sum_{j\in\I_0}\bI(W_j\leq -t)}{{1\vee\sum_j\bI(W_j\leq -t)}}(D_1+D_2+D_3+D_4)
\end{align*}

By Lemmas \ref{klem1}-\ref{klem4}, we get
\begin{align*}
D_1&=O_p\left(\sqrt{\xi_{n,p}/n}+\zeta_p\right)-\frac{2t^3\bar{\kappa}}{9n},\\
D_2&=O_p\left(\sqrt{\xi_{n,p}/n}+\zeta_p+(\log p)^{-1}\right)-\frac{5t^3\bar{\kappa}}{18n}.
\end{align*}

Next, we deal with $D_3$, and the $D_4$ can be bounded in a similar way to $D_3$.  Let $\mathcal{C}_p$ denote the event that $\left\{\max_{j\in \I_1}|s_{kj}^2/\sigma^2_j-1|<\delta_n,k=1,2\right\}$ and $\delta_n=c\sqrt{n^{-1}\log p}$ for some sufficiently large $c>0$.
By Assumption \ref{moment}-(ii), we have
\begin{align*}
\Pr\left\{\max_{j\in\I_1}|s_{kj}^2/\sigma^2_j-1|\geq\delta_n\right\}=O(\sqrt{\xi_{n,p}/n}),
\end{align*}
and thus ${\Pr(\mathcal{C}_p^c)}=O(\sqrt{\xi_{n,p}/n})$ holds.

Conditional on the event $\mathcal{C}_p$, we have
\begin{align*}
|D_3|&=\frac{\left|\sum_{j\in\I_1}\bI(W_j\geq t )-\sum_{j\in\I_1}\bI(\tW_j\geq t)\right|}{\sum_{j\in\I_0}\bI(W_j\leq -t )}\\
&\leq \frac{\sum_{j\in\I_1}\left|\bI(W_j\geq t )-\bI(\tW_j\geq t)\right|}{\sum_{j\in\I_0}\bI(W_j\leq -t )}\\
&=\frac{\sum_{j\in\I_1}\bI(t\leq W_j\leq t(1+3\delta_n))}{{\sum_{j\in\I_0}\bI(W_j\leq -t )}}+\frac{\sum_{j\in\I_1}\bI(t(1-3\delta_n)\leq W_j\leq t)}{{\sum_{j\in\I_0}\bI(W_j\leq -t )}}\\
&\equiv D_{31}+D_{32}.
\end{align*}
We only need to deal with $D_{31}$ and $D_{32}$ follows similarly. By the Markov's inequality, we get
\begin{align*}
\Pr(D_{31}>\epsilon)\leq& p_0^{-1}\frac{\sum_{j\in \I_1}\Pr(t\leq W_j\leq t(1+3\delta_n))}{\epsilon G(t)(1+\zeta_p)}\\
=&p_0^{-1}\frac{\sum_{j\in \I_1\setminus\mathcal{S}_p}\Pr(t\leq W_j\leq t(1+3\delta_n))+\sum_{j\in\mathcal{S}_p}\Pr(t\leq W_j\leq t(1+3\delta_n))}{\epsilon G(t)(1+\zeta_p)}\\
\leq& p_0^{-1}\frac{\sum_{j\in \I_1\setminus\mathcal{S}_p}\Pr(t\leq W_j\leq t(1+3\delta_n))+C\beta_p p^{-1}}{\epsilon G(t)(1+\zeta_p)}\\
\leq& p_0^{-1}\frac{\sum_{j\in \I_1\setminus\mathcal{S}_p}\Pr(t\leq W_j\leq t(1+3\delta_n))+C\beta_p p^{-1}}{\epsilon t^{-1}\exp(-t)(1+\zeta_p)},\\
\leq& C\frac{(p_1-\beta_p)t^2\exp(t)\delta_n+\beta_pp^{-1}t\exp(t)}{p_0\epsilon(1+\zeta_p)}\\
\leq& C\frac{(p_1/\beta_p-1)\log^2(p/\beta_p)\delta_n}{\epsilon(1+\zeta_p)}+O(p^{-1}\log p)=O(\delta_n),
\end{align*}
where $\mathcal{S}_p=\left\{j: |\mu_j|/\sigma_j\geq 3\sqrt{\log p/n}\right\}$, and $C$ is some positive constant. The second equality is due to Lemma \ref{klm2}, the fourth inequality comes from the fact that $t\lesssim \log\{p/(\beta_p\alpha)\}$ and the last inequality uses the condition $(p_1/\beta_p-1)\log^2(p/\beta_p)=O(1)$ .

Finally, collecting all the terms of $D_k,k=1,\ldots,4$, we conclude that for any $G^{-1}(\alpha p_1/\{p_0(1-\alpha)\})\lesssim t\lesssim \log\{p/(\beta_p\alpha)\}$,
\[
\theta(t)=\frac{\sum_{j\in\I_0}\bI(W_j\leq -t)}{{1\vee\sum_j\bI(W_j\leq -t)}}\left[-\frac{t^3\bar{\kappa}}{2n}+O_p\left(\sqrt{\xi_{n,p}/n}+\zeta_p\right)\right].
\]
By using similar arguments given in the Supplemental Material we can show that $\frac{\sum_{j\in\I_0}\bI(W_j\leq -t)}{{1\vee\sum_j\bI(W_j\leq -t)}}=1+O_p(p^{\eta-1})$. Thus, we conclude that $\theta(t)=-\frac{t^3\bar{\kappa}}{2n}\left\{1+O_p(p^{\eta-1})\right\}+O_p\left(\sqrt{\xi_{n,p}/n}+\zeta_p\right)$.


Accordingly, we obtain
\begin{align*}
&\FDP_W(\wht)=\frac{\sum_{j\in\I_0}\bI(W_j\geq \wht)}{1\vee\sum_j\bI(W_j\geq \wht)}\\
&=\frac{\sum_{j}\bI\left(W_j\leq- \wht \right)}{1\vee\sum_j\bI(W_j\geq \wht)}\times\frac{\sum_{j\in\I_0}\bI(W_j\geq \wht)}{\sum_{j}\bI\left(W_j\leq- \wht \right)}\\
&=\frac{\sum_{j}\bI\left(W_j\leq- \wht \right)}{1\vee\sum_j\bI(W_j\geq \wht)}\times\frac{\sum_{j\in\I_0}\bI(W_j\geq \wht)}{\sum_{j\in\I_0}\bI\left(W_j\leq- \wht \right)}\times\frac{\sum_{j\in\I_0}\bI(W_j\leq -\wht)}{\sum_{j}\bI\left(W_j\leq- \wht \right)}\\
&=\frac{\sum_{j}\bI\left(W_j\leq- \wht \right)}{1\vee\sum_j\bI(W_j\geq \wht)}\left\{1+\frac{2t^3\bar{\kappa}}{9n}+O_p\left(\sqrt{\xi_{n,p}/n}+\zeta_p\right)\right\}\left\{1+O_p(p^{\eta-1})\right\}\\\
&=\frac{\sum_{j}\bI\left(W_j\leq- \wht \right)}{1\vee\sum_j\bI(W_j\geq \wht)}\left\{1-\frac{4}{9}\theta(\wht)\right\}\left\{1+O_p\left(\sqrt{\xi_{n,p}/n}+\zeta_p\right)\right\}\\
&\leq \alpha+O_p\left(\sqrt{\xi_{n,p}/n}+\zeta_p\right).
\end{align*}

The proof is completed. \hfill$\Box$

\bibliographystyle{asa}
\bibliography{Mirror-ref}

\newpage

{\Large{\bf Supplementary Material for
``A New Procedure for Controlling False Discovery Rate in Large-Scale $t$-tests"}}

%

\baselineskip 20pt

\def\thelemma{S.\arabic{lemma}}
\def\thepro{S.\arabic{pro}}
\def\theequation{S.\arabic{equation}}
\def\thetable{S\arabic{table}}
\def\thefigure{S\arabic{figure}}
\setcounter{lemma}{0}
\setcounter{figure}{0}
\setcounter{table}{0}

This supplementary material contains the proofs of some technical lemmas and corollaries, and additional simulation results.

\subsection*{Additional Lemmas}

The first one is the large deviation result for the Student's $t$ statistic $T$. See also \cite{wang2005}.

\begin{lemma}\label{ldp}[\cite{delaigle2011robustness}]
Let $B>1$ denote a constant. Then,
\[
\frac{\Pr(T>x)}{1-\Phi(x)}=\exp\left\{-\frac{x^3\kappa}{3\sqrt{n}}\right\}\left[1+\theta(n,x)\left\{(1+|x|)n^{-1/2}+(1+|x|)^4n^{-1}\right\}\right]
\]
as $n\to\infty$, where the function $\theta$ is bounded in absolute value by a finite, positive constant $C_1(B)$ (depending only on $B$), uniformly in all distributions of $X$ for which $\E|X|^4\leq B$, $\E(X^2)=1$ and $\E(X) = 0$, and uniformly in $x$ satisfying $0\leq x \leq Bn^{1/4}$.
\end{lemma}

The second one is a standard large deviation result for the mean; See Theorem VIII-4 in \cite{petrov2012sums}.
\begin{lemma}[Large deviation for the mean]\label{ldm}
Suppose that $X_1,\ldots,X_n$ are i.i.d random variables with mean zero and variance $\sigma^2$, satisfying Assumption \ref{moment}-(ii). Then for any $0\leq x\leq cn^{1/6}$ and $c>0$,
\[
\frac{\Pr(\sqrt{n}\bar{X}/\sigma>x)}{1-\Phi(x)}=\exp\left\{\frac{x^3\kappa}{6\sqrt{n}}\right\}\left\{1+o(1)\right\}.
\]
\end{lemma}

The third lemma is a large deviation result for $T_{j}/(s_j/\sigma_j)$.
\begin{lemma}\label{ldtw}
Suppose Assumptions \ref{moment}-(ii) hold. Then for $x\to\infty$ and $x=o(n^{1/6})$,
\[
\frac{\Pr(T/(s/\sigma)>x)}{1-\Phi(x)}=\exp\left\{-\frac{5x^3\kappa}{6\sqrt{n}} \right\}\left\{1+O\left(x^{-2}+x/\sqrt{n}\right)\right\}.
\]
\end{lemma}

\proof Without loss of generality, we assume that $\sigma^2=1$.  First of all, we deal with $\Pr(\sqrt{n}\bar{X}/m_2>x)$, where $m_2=n^{-1}\sum_{i=1}^nX_i^2$. Observe
\begin{align*}
\Pr(\sqrt{n}\bar{X}/m_2>x)=\Pr\left(n^{-1/2}\sum_{i=1}^n Y_i\geq x\right),
\end{align*}
where $Y_i=X_i-c_{n,x}(X_i^2-1)$ and $c_{n,x}=n^{-1/2}x$. Simple calculation yields
\begin{align*}
\var(Y_i)&=(1-2c_{n,x}\kappa)+O(c^2_n),\\
E(Y_i^3)&=\kappa-3c_{n,x}\E X_i^4+3c_{n,x}+O(c^2_n).
\end{align*}

By Lemma \ref{ldm}, we have
\begin{align*}
\frac{\Pr\left(\frac{n^{-1/2}\sum_{i=1}^n Y_i}{\sqrt{\var(Y)}}\geq \frac{x}{\sqrt{\var(Y)}}\right)}{\tilde{\Phi}(x/\sqrt{\var(Y)})}&=\exp\left\{\frac{x^3\kappa_{Y}}{6\sqrt{n}\var^{3/2}{Y}}\right\}\left\{1+O(x/\sqrt{n})\right\}\\
&=\exp\left\{\frac{x^3\kappa}{6\sqrt{n}}\right\}\left\{1+O(x/\sqrt{n})\right\},
\end{align*}
where $\kappa_Y=E(Y_i^3)/\var^{3/2}(Y_i)$.

By using the fact that for $x\to\infty$ and $x=o(a_n^{-1/2})$,
\[
\tilde{\Phi}(x(1+a_n))=\tilde{\Phi}(x)\exp(-x^2{a_n})\left\{1+O(x^{-2})\right\},
\]
we have
\begin{align*}
{\tilde{\Phi}(x/\sqrt{\var(Y)})}={\tilde{\Phi}(x)}\exp\left\{-x^2c_{n,x}\kappa \right\}\left\{1+O(x^{-2})\right\}.
\end{align*}
Thus,
\begin{align*}
\Pr(\sqrt{n}\bar{X}/m_2>x)&={\tilde{\Phi}(x)}\exp\left\{-x^2c_{n,x}\kappa+\frac{x^3\kappa}{6\sqrt{n}} \right\}\left\{1+O\left(x^{-2}+x/\sqrt{n}\right)\right\}\\
&={\tilde{\Phi}(x)}\exp\left\{-\frac{5x^3\kappa}{6\sqrt{n}} \right\}\left\{1+O\left(x^{-2}+x/\sqrt{n}\right)\right\}.
\end{align*}

Finally, we show that ${\Pr(T/s>x)}$ and $\Pr(\sqrt{n}\bar{X}/m_2>x)$ are close enough. Note that
\begin{align*}
\Pr(T/s>x)&=\Pr(\sqrt{n}\bar{X}/(m_2-\bar{X}^2)>x)\\
&=\Pr(x\bar{X}^2+\sqrt{n}\bar{X}-m_2x>0)\\
&=\Pr\left(\frac{\sqrt{n}\bar{X}}{m_2}>\frac{-{n}+\sqrt{n^2+4nx^2m_2}}{2xm_2}\right)+\Pr\left(\frac{\sqrt{n}\bar{X}}{m_2}<\frac{-{n}-\sqrt{n^2+4nx^2m_2}}{2xm_2}\right)\\
&=\Pr\left(\frac{\sqrt{n}\bar{X}}{m_2}>x\left\{1+O(x^2n^{-1})\right\}\right)+o(x/\sqrt{n})\tilde{\Phi}(x)\\
&=\tilde{\Phi}(x\left\{1+O(x^2n^{-1})\right\})\exp\left\{-\frac{5x^3\kappa}{6\sqrt{n}} \right\}\left\{1+O\left(x^{-2}+x/\sqrt{n}\right)\right\}\\
&=\tilde{\Phi}(x)\exp\left\{-\frac{5x^3\kappa}{6\sqrt{n}} \right\}\left\{1+O\left(x^{-2}+x/\sqrt{n}\right)\right\}.
\end{align*}
\hfill$\Box$

\subsection*{Proof of Lemmas and Propositions}

\noindent{\bf Proof of Lemma \ref{klem1}.} Recalling $p=\exp\{o(n^{1/3})\}$, we have
$t=o(n^{1/3})$ and $a_{nt}=o(n^{1/6})$. Let $\A_p=\{v: t/a_{nt}<|v|<a_{nt}\}$. Then,
\begin{align*}
&\frac{\Pr(T_1T_2>t)}{\Pr(T_1T_2<-t)}-1\\
&=\frac{\Pr(T_1T_2>t,\A_p)-\Pr(T_1T_2<-t,\A_p)}{\Pr(T_1T_2<-t)}+\frac{\Pr(T_1T_2>t,\A^c_p)-\Pr(T_1T_2<-t,\A^c_p)}{\Pr(T_1T_2<-t)}\\
&\equiv C_1+C_2.
\end{align*}

Firstly, for the term $C_2$,
\begin{eqnarray*}
C_2&=&\left[\frac{\Pr(T_1T_2>t,\A^c_p)}{P(Z_1Z_2>t)}-\frac{\Pr(T_1T_2<-t,\A^c_p)}{\Pr(Z_1Z_2<-t)}\right]\frac{\Pr(Z_1Z_2<-t)}{\Pr(T_1T_2<-t)},
\end{eqnarray*}
where $Z_1$ and $Z_2$ are two independent $N(0,1)$ variables. From the proof given later, it can be easily see that ${\Pr(Z_1Z_2<-t)}/{\Pr(T_1T_2<-t)}\rightarrow 1$ uniformly in $0\leq t\leq C\log p$. Thus, in what follows we mainly focus on the rate of ${\Pr(T_1T_2>t,\A^c_p)}/{\Pr(Z_1Z_2>t)}$. The other term $\Pr(T_1T_2<-t,\A^c_p)$ can be handled similarly.

Note that
$\Pr(T_1T_2>t, \A^c_p)\leq 2\Pr(|T_1|>a_{nt})$, $\Pr(Z_1Z_2>t)\geq \left\{\Pr(Z_1>\sqrt t)\right\}^2$. By the inequality
\[
\frac{x}{x^2+1}\phi(x)<\tilde{\Phi}(x)<\phi(x)/x, \ \ \mbox{for all} \  x
\]
and the large deviation formula for the $t$-statistic (Lemma \ref{ldp}), we obtain that
\begin{align*}
\frac{\Pr(T_1T_2>t, \A^c_p)}{\Pr(Z_1Z_2>t)}&\leq \frac{2\Pr(|T_1|>a_{nt})}{\left\{\Pr(Z_1>\sqrt t)\right\}^2}\\
&\leq c \exp\left\{-\frac{1}{2}(a_{nt}^2-t-t)\right\}\frac{t}{a_{nt}}\\
&=O\left(\sqrt{t/n}\right),
\end{align*}
thus we claim that $C_2=O(\sqrt{t/n})$.

Next, we deal with the main term $C_1$. Denote $\A_p^{+}=\{v: t/a_{nt}<v<a_{nt}\}$.  Observe
\begin{align*}
&\Pr(T_1>t/T_2, \A_p^{+})\\
=&-\int_{\A^{+}_p} \tilde{\Phi}({t}/{v})\exp\left\{-\frac{t^3\kappa}{3v^3\sqrt n}\right\}\left\{1+O({t}/{(v\sqrt{n})})\right\}d\tilde{\Phi}(v)\exp\left\{-\frac{v^3\kappa}{3\sqrt n}\right\}\left\{1+O(v/\sqrt{n})\right\}\\
=&\int_{\A^{+}_p} \tilde{\Phi}({t}/{v})\phi(v)\exp\left\{-\frac{\{(t/v)^3+v^3\}\kappa}{3\sqrt n}\right\}dv\left\{1+O(a_{nt}/\sqrt{n})\right\}\\
&+\int_{\A^{+}_p} \tilde{\Phi}({t}/{v})\tilde{\Phi}(v)\exp\left\{-\frac{\{(t/v)^3+v^3\}\kappa}{3\sqrt n}\right\}\frac{v^2\kappa}{\sqrt{n}}dv\left\{1+O(a_{nt}/\sqrt{n})\right\}\\
\equiv &R_1+R_2,
\end{align*}
where we use Lemma \ref{ldp} again. Note that
\begin{align*}
R_2&\leq \int_{\A^{+}_p}\tilde{\Phi}({t}/{v})\phi(v)\exp\left\{-\frac{\{(t/v)^3+v^3\}\kappa}{3\sqrt n}\right\}\frac{v\kappa}{\sqrt{n}}dv\left\{1+O(a_{nt}/\sqrt{n})\right\},
\end{align*}
and accordingly $R_2=R_1O(a_{nt}/\sqrt{n})$. Hence,
\[
\Pr(T_1T_2>t, \A_p^{+})=\int_{\A^{+}_p} \tilde{\Phi}({t}/{v})\phi(v)\exp\left\{-\frac{\{(t/v)^3+v^3\}\kappa}{3\sqrt n}\right\}dv\left\{1+O(a_{nt}/\sqrt{n})\right\}.
\]

Similarly,
\begin{align*}
\Pr(T_1T_2>t, \A_p^{-})&=\int_{\A^{+}_p} \tilde{\Phi}({t}/{v})\phi(v)\exp\left\{\frac{\{(t/v)^3+v^3\}\kappa}{3\sqrt n}\right\}dv\left\{1+O(a_{nt}/\sqrt{n})\right\},\\
\Pr(T_1T_2<-t, \A_p^{+})&=\int_{\A^{+}_p} \tilde{\Phi}({t}/{v})\phi(v)\exp\left\{\frac{\{(t/v)^3-v^3\}\kappa}{3\sqrt n}\right\}dv\left\{1+O(a_{nt}/\sqrt{n})\right\},\\
\Pr(T_1T_2<-t, \A_p^{-})&=\int_{\A^{+}_p} \tilde{\Phi}({t}/{v})\phi(v)\exp\left\{\frac{\{-(t/v)^3+v^3\}\kappa}{3\sqrt n}\right\}dv\left\{1+O(a_{nt}/\sqrt{n})\right\}.
\end{align*}
The Taylor's expansion yields that
\begin{align*}
\Pr(T_1T_2>t)-\Pr(T_1T_2<-t)=\frac{4t^3\kappa^2}{9n}\int_{\A^{+}_p} \tilde{\Phi}({t}/{v})\phi(v)dv\left\{1+O(a_{nt}/\sqrt{n})\right\}.
\end{align*}
Consequently, we easily get that
\begin{eqnarray*}
C_1=\frac{2t^3\kappa^2}{9n}+O(a_{nt}/\sqrt{n}),
\end{eqnarray*}
from which we obtain the assertion. \hfill$\Box$

\noindent{\bf Proof of Lemma \ref{klmwt}.} The proof is similar to Lemma \ref{klem1} but
using the large deviation result for $T/(s/\sigma)$ obtained in Lemma \ref{ldtw} and for
the mean. \hfill$\Box$

\noindent{\bf Proof of Lemma \ref{klm2}.} We prove (\ref{mco1}) and the proof of
(\ref{mco2}) follows similarly. Here $W_j=T_{1j}T_{2j}$. Clearly, $G(t)$ is a deceasing
and continuous function. Let $z_0<z_1<\cdots<z_{d_p}\leq 1$ and $t_i=G^{-1}(z_i)$, where
$z_0=b_p/p, z_i=b_p/p+c_p\exp(i^{\delta})/p, d_p=[\{\log((p-b_p)/c_p)\}^{1/\delta}]$ with
$c_p/b_p\to 0$ and $0<\delta<1$. Note that $G(t_i)/G(t_{i+1})=1+o(1)$ uniformly in $i$.
Thus, it is enough to obtain the convergence rate of
\[
D_p=\sup_{0\leq i\leq d_p}\left|\frac{\sum_{j\in \I_0}\left\{\bI(W_j>t_i)-\Pr(W_j>t_i)\right\}}{p_0G(t_i)}\right|.
\]

Define $\mathcal{S}_j=\{k\in\I_0:\,\, X_k\,\,\mbox{is dependent with}\,\, X_j\}$ and further
\[
D(t)=\E\left[\sum_{j\in \I_0}\left\{\bI(W_j>t)-\Pr(W_j>t)\right\}^2\right].
\]

It is noted that
\begin{align*}
D(t)&=\sum_{j\in \I_0}\sum_{k\in \mathcal{S}_j}\E\left[\left\{\bI(W_j>t)-\Pr(W_j>t)\right\}\left\{\bI(W_k>t)-\Pr(W_k>t)\right\}\right]\\
&\leq  r_pp_0G(t).
\end{align*}

We then obtain
\begin{align*}
\Pr(D_p\geq \epsilon)&\leq \sum_{i=0}^{d_p}
\Pr\left(\left|\frac{\sum_{j\in \I_0} [\bI(W_j>t_i)-\Pr(W_j>t_i)]}{p_0G(t_i)}\right|\geq \epsilon\right)\\
&\leq  \frac{1}{\epsilon^2}\sum_{i=0}^{d_p}\frac{1}{p^2_0G^2(t_i)}D(t_i)\\
&\leq  \frac{r_p}{\epsilon^2}\sum_{i=0}^{d_p}\frac{1}{p_0G(t_i)}.
\end{align*}

Moreover, observe that
\begin{eqnarray*}
&&\sum_{i=0}^{d_p}\frac{1}{p_0G(t_i)}=\frac{p}{p_0}\left(\frac{1}{b_p}
+\sum_{i=1}^{d_p}\frac{1}{b_p+c_pe^{i^{\delta}}}\right)\\
&\leq & c\left(\frac{1}{b_p}+c_p^{-1}\sum_{i=1}^{d_p}\frac{1}{1+e^{i^{\delta}}}\right)
\leq c/c_p\{1+o(1)\}.
\end{eqnarray*}
Because $c_p$ can be made arbitrarily large as long as $c_p/b_p\to 0$, we have $D_p=O_p(\sqrt{r_p/b_p})$.

%

\noindent{\bf Proof of Lemma \ref{klem4}.} By using the same arguments given in the proof
of Lemma \ref{klm2}, this lemma can be proved easily and thus omitted.   \hfill$\Box$
\subsubsection*{Proof of Proposition \ref{pro2}}

We prove this proposition for $L_{+}$. The result for $L$ can be obtained similarly. Fix $\epsilon>0$ and for any threshold $t>0$, define
\[
R_{\epsilon}(t)=\frac{\sum_{j\in\I_0}\bI\left(W_j\geq t,\Delta_j\leq\epsilon\right)}{1+\sum_{j\in\I_0}\bI\left(W_j\leq -t\right)}.
\]
Consider the event that $\mathcal{A}=\{\Delta\equiv\max_{j\in\I_0}\Delta_j\leq\epsilon\}$. Furthermore, for a threshold rule $L=T({\bf W})$ mapping statistics ${\bf W}$ to a threshold $L\geq 0$, for each index $j=1,\ldots,p$, we define
\[
L_j=T\left(W_1,\ldots,W_{j-1},|W_j|,W_{j+1},\ldots,W_{p}\right)\geq 0
\]
i.e. the threshold that we would obtain if $\sgn(W_j)$ were set to 1.

Then for the RESS method with the threshold $L_{+}$, we can write
\begin{align*}
\frac{\sum_{j\in\I_0}\bI\left(W_j\geq L_{+},\Delta_j\leq\epsilon\right)}{1\vee\sum_j\bI(W_j\geq L_{+})}&=\frac{1+\sum_{j}\bI\left(W_j\leq-L_{+}\right)}{1\vee\sum_j\bI(W_j\geq L_{+})}\times\frac{\sum_{j\in\I_0}\bI\left(W_j\geq L_{+},\Delta_j\leq\epsilon\right)}{1+\sum_{j}\bI\left(W_j\leq-L_{+}\right)}\\
&\leq \alpha\times R_{\epsilon}(L_{+}).
\end{align*}
It is crucial to get an upper bound for $\E\{R_{\epsilon}(L_{+})\}$. We have
\begin{align}
\E\{R_{\epsilon}(L)\}&=\E\left[\frac{\sum_{j\in\I_0}\bI\left(W_j\geq L,\Delta_j\leq\epsilon\right)}{1+\sum_{j\in\I_0}\bI\left(W_j\leq -L_{+}\right)}\right]\nonumber\\
&=\sum_{j\in\I_0}\E\left\{\frac{\bI\left(W_j\geq L_j,\Delta_j\leq\epsilon\right)}{1+\sum_{k\in\I_0,k\neq j}\bI\left(W_k\leq -L_j\right)}\right\} \nonumber\\
&=\sum_{j\in\I_0}\E\left[\E\left\{\frac{\bI\left(W_j\geq L_j,\Delta_j\leq\epsilon\right)}{1+\sum_{k\in\I_0,k\neq j}\bI\left(W_k\leq -L_j\right)}\mid |W_j|\right\}\right]\nonumber \\
&=\sum_{j\in\I_0}\E\left\{\frac{\Pr\left(W_j>0\mid |W_j|\right)\bI\left(|W_j|\geq L_j,\Delta_j\leq\epsilon\right)}{1+\sum_{k\in\I_0,k\neq j}\bI\left(W_k\leq -L_j\right)}\right\},\label{kf}
\end{align}
where the last step holds since the only unknown is the sign of $W_j$ after conditioning on $|W_j|$. By definition of $\Delta_j$, we have
$\Pr\left(W_j>0\mid |W_j|\right)\leq 1/2+\Delta_j$.

Hence,
\begin{align*}
\E&\{R_{\epsilon}(L_{+})\}\\
&\leq \sum_{j\in\I_0}\E\left\{\frac{(\frac12+\Delta_j)\bI\left(|W_j|\geq L_j,\Delta_j\leq\epsilon\right)}{1+\sum_{k\in\I_0,k\neq j}\bI\left(W_k\leq -L_j\right)}\right\}\\
&\leq(\frac12+\epsilon)\left[\sum_{j\in\I_0}\E\left\{\frac{\bI\left(W_j\geq L_j,\Delta_j\leq\epsilon\right)}{1+\sum_{k\in\I_0,k\neq j}\bI\left(W_k\leq -L_j\right)}\right\}+\sum_{j\in\I_0}\E\left\{\frac{\bI\left(W_j\leq -L_j\right)}{1+\sum_{k\in\I_0,k\neq j}\bI\left(W_k\leq -L_j\right)}\right\}\right]\\
&=(\frac12+\epsilon)\left[\E\{R_{\epsilon}(L_{+})\}+\sum_{j\in\I_0}\E\left\{\frac{\bI\left(W_j\leq -L_j\right)}{1+\sum_{k\in\I_0,k\neq j}\bI\left(W_k\leq -L_j\right)}\right\}\right].
\end{align*}
Finally, the sum in the last expression can be simplified as: if for all null $j$, $W_j>-L_j$, then the sum is equal to zero, while otherwise,
\begin{align*}
\sum_{j\in\I_0}\E\left\{\frac{\bI\left(W_j\leq -L_j\right)}{1+\sum_{k\in\I_0,k\neq j}\bI\left(W_k\leq -L_j\right)}\right\}=\sum_{j\in\I_0}\E\left\{\frac{\bI\left(W_j\leq -L_j\right)}{1+\sum_{k\in\I_0,k\neq j}\bI\left(W_k\leq -L_k\right)}\right\}=1,
\end{align*}
where the first step comes from the fact: for any $j,k$, if $W_j\leq-\min(L_j,L_k)$ and $W_k\leq-\min(L_j,L_k)$, then $L_j=L_k$. See \cite{barber2018robust} for a proof.

Accordingly, we have
\[
\E\{R_{\epsilon}(L_{+})\}\leq \frac{1/2+\epsilon}{1/2-\epsilon}\leq 1+4\epsilon.
\]
 Consequently, the assertion of this proposition holds. \hfill$\Box$

\noindent{\bf Proof of Proposition \ref{pro3}.} The proof follows similarly to that of
Proposition \ref{pro2} but uses $\Pr\left(W_j>0\mid |W_j|, {\bf W}_{-j}\right)$ to
replace $\Pr\left(W_j>0\mid |W_j|\right)$ in (\ref{kf}).  \hfill$\Box$

\noindent{\bf Discussion on $R(t)$.} Recall the definition of $\mathcal{M}=\left\{j:
|\mu_j|/\sigma_j\geq \sqrt{(2+c)\log p/n}\right\}$ for any $c>0$. We shall show that
\[
\sup_{0\leq t\leq 2\log p}\left|\frac{\sum_{j\in\I_1}\bI\left(W_j\leq -t\right)}{\sum_{j\in\I_0}\bI\left(W_j\leq -t\right)}\right|\cp 0.
\]

Observe that for $0\leq t\leq 2\log p$,
\begin{align*}
&\frac{\sum_{j\in\mathcal{M}}\Pr\left(W_j\leq -t\right)}{\sum_{j\in\I_0}\Pr\left(W_j\leq -t\right)}\\
&\leq \frac{2\sum_{j\in\mathcal{M}}\Pr(T_{1j}\leq -\sqrt{t})}{p_0\exp(-t)/t\left\{1+o(1)\right\}}\\
&\leq\frac{2|\mathcal{M}|\max_{j\in\I_0}\Pr(T_{1j}\leq -\sqrt{t}-\sqrt{(1+c)\log p})}{p_0\exp(-t)/t\left\{1+o(1)\right\}}=o(1).
\end{align*}
On the other hand, $\sum_{j\in \I_1\backslash\mathcal{M}}\Pr\left(W_j\leq -t\right)\leq (p_1-|\mathcal{M}|)\Pr_{j\in\I_0}\left(W_j\leq -t\right)$ and thus the assertion holds due to $(p_1-|\mathcal{M}|)/p_0\to 0$.

\begin{table}    \tabcolsep 9pt
\caption{\small Comparison results of FDR and TPR when the signals are $\mu_j=\delta_j\sqrt{\log p/n_t}$ with $\delta_j\sim {\rm Unif}(-1.5,-1)$ under $p=5000$, $\rho=0$, $\vartheta=0.05$ and the target FDR $\alpha=0.2$.} \label{tab:one}\vspace{0.3cm}
\linespread{1.2}
\centering
{\small\renewcommand{\arraystretch}{1}
\begin{tabular}{cccccccccccccccccc}
\toprule
&&  \multicolumn{2}{c}{$n_t=50$} && \multicolumn{2}{c}{$n_t=100$} \\
& Method & \tc{FDR(\%)} & \tc{TPR(\%)} && \tc{FDR(\%)} & \tc{TPR(\%)} \\
\hline
$t(5)$  &   RESS  &   $20.1_{\tiny{(6.0)}}$    &   $64.0_{\tiny{(6.0)}}$ && $21.0_{\tiny{(6.2)}}$    &   $65.2_{\tiny{(6.4)}}$\\
&   BH  &   $17.3_{\tiny{(2.8)}}$    &   $63.6_{\tiny{(3.6)}}$    &&   $18.2_{\tiny{(2.7)}}$    &   $65.4_{\tiny{(3.5)}}$\\
\\
exp(1)  &   RESS  &   $32.9_{\tiny{(4.0)}}$    &   $74.1_{\tiny{(3.5)}}$    &&  $27.5_{\tiny{(3.4)}}$    &   $77.4_{\tiny{(3.3)}}$\\
&   BH  &   $49.1_{\tiny{(2.1)}}$    &    $85.6_{\tiny{(2.4)}}$ && $41.1_{\tiny{(2.9)}}$    &   $86.7_{\tiny{(2.5)}}$\\
\\
Mixed   &   RESS  &  $26.0_{\tiny{(4.5)}}$    &   $74.3_{\tiny{(4.3)}}$    &&  $23.5_{\tiny{(4.6)}}$    &   $76.2_{\tiny{(4.0)}}$\\
&   BH  &   $31.6_{\tiny{(2.6)}}$    &   $80.3_{\tiny{(2.6)}}$    &&  $27.6_{\tiny{(3.1)}}$    &   $80.7_{\tiny{(2.5)}}$\\
\bottomrule
\end{tabular}}
\end{table}

\subsection*{Additional simulation results}

Table S1 reports a brief comparison between this RESS and BH procedures when $\mu_j<0$ for all $j\in\mathcal{I}_1$. We see that the RESS significantly improves the accuracy of FDR control over the BH to certain degree, while maintains high power in most cases.

\begin{figure}[ht]\centering
\includegraphics[width=1.0\textwidth]{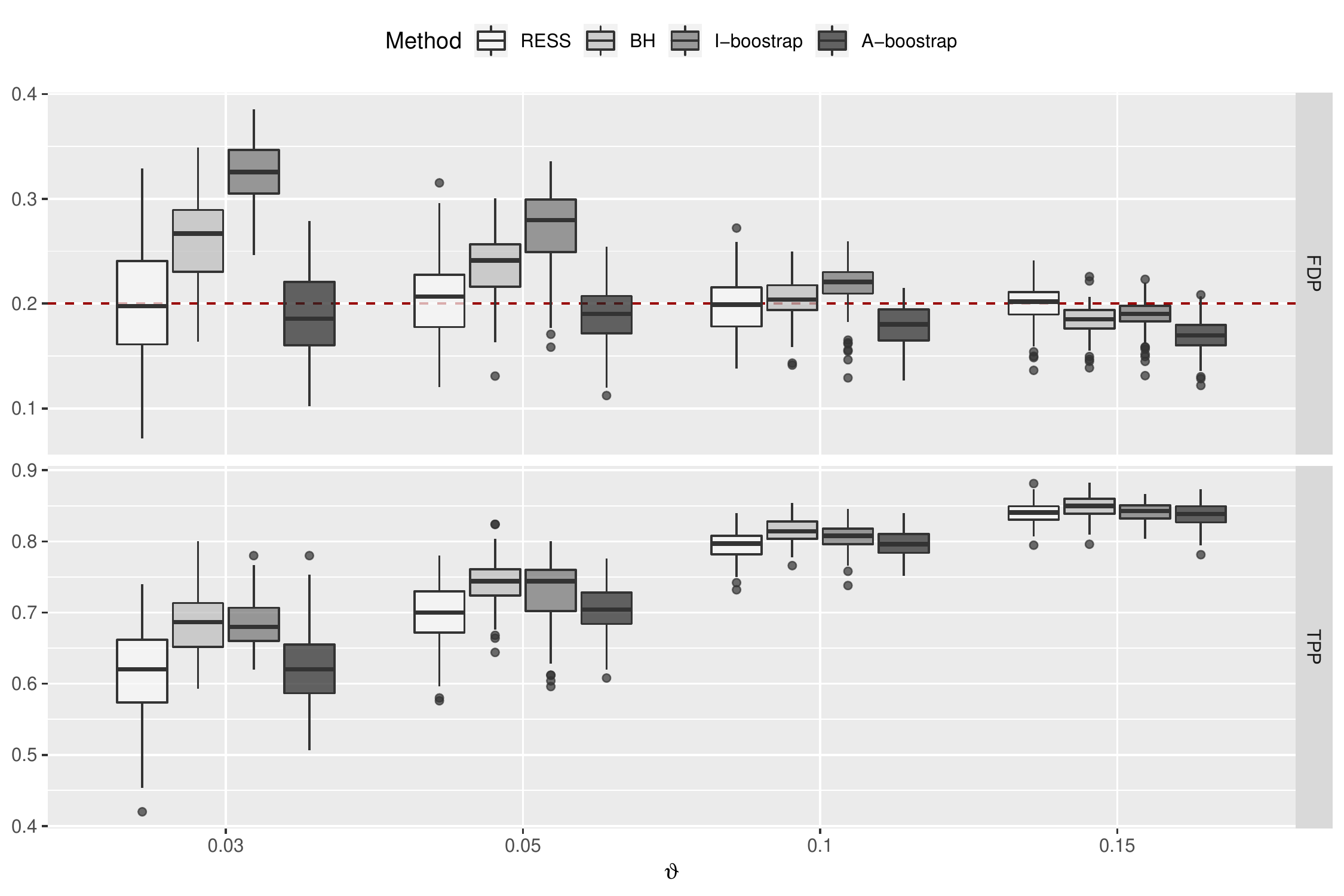}\vspace{-0.2cm}
\caption {\small Box-plots of false discovery proportion (FDP) and
true positive proportion (TPP) when errors are distributed from the
mixed distribution under $n_t=100$, $p=5000$, and $\rho=0$;
The red dashed lines indicate the target FDR level.}\label{Fig:ratio2}\vspace{-0.2cm}
\end{figure}

Figure \ref{Fig:ratio2} is a counterpart of Figure \ref{Fig:ratio} when the errors are distributed from the mixed distribution.

\end{document}